\documentclass[11pt,a4paper]{article}
 
\usepackage[T1]{fontenc}
\usepackage[utf8]{inputenc}
\usepackage[british]{babel}
\usepackage{lmodern}
\usepackage{amsmath,amssymb,mathtools,amsthm}
\usepackage{microtype}
\usepackage[a4paper,margin=3cm]{geometry}
 
\usepackage{graphicx}
\usepackage{tikz}
\usetikzlibrary{arrows.meta,positioning,fit,calc}
\usepackage{pgfplots}
\pgfplotsset{compat=1.18}
\pgfplotscreateplotcyclelist{greys}{%
  {fill=black!85,draw=black},{fill=black!60,draw=black},{fill=black!38,draw=black},%
  {fill=black!18,draw=black},{fill=white,draw=black}}
\usepackage{booktabs}
\usepackage{tabularx}
\usepackage{longtable}
\usepackage{array}
\usepackage[table]{xcolor}
\usepackage{threeparttable}
\usepackage{caption}
\usepackage{subcaption}
\usepackage{float}
\usepackage{pdflscape}
\usepackage{enumitem}
\usepackage{siunitx}
\usepackage{csquotes}
\usepackage[
  backend=biber,
  style=authoryear,
  maxcitenames=2,
  maxbibnames=99,
  giveninits=true,
  useprefix=true,
  uniquename=init,
  doi=true,
  url=true,
  isbn=false
]{biblatex}
\usepackage{algorithm}
\usepackage{algpseudocode}
 
\usepackage{hyperref}
\usepackage[nameinlink,capitalise,noabbrev]{cleveref}
\hypersetup{
  hidelinks,
  pdftitle={A Strategic Freight Market Model for Liner Shipping},
  pdfauthor={Jonas Westin},
  pdfkeywords={freight transport, liner shipping, Bertrand competition, discrete choice, nested logit, market entry, Ramsey pricing}
}
 
\theoremstyle{plain}
\newtheorem{property}{Property}[section]
\crefname{property}{Property}{Properties}
 
\numberwithin{equation}{section}
\renewcommand{\arraystretch}{1.1}

\newcommand{\SO}{\mathrm{SO}}

\newcommand{\fmBaseCompPS}{\num{190827}}
\newcommand{\fmBaseCompShare}{58.7}
\newcommand{\fmBaseCompPrice}{896}
\newcommand{\fmBaseCompN}{4}

\newcommand{\fmBaseCompLF}{0.476}
\newcommand{\fmBaseMonoW}{\num{2725430}}

\newcommand{\fmBaseMonoPS}{\num{336262}}
\newcommand{\fmBaseMonoShare}{45.5}
\newcommand{\fmBaseMonoPrice}{1019}
\newcommand{\fmBaseMonoN}{3}
\newcommand{\fmBaseMonodW}{\num{-208345}}
\newcommand{\fmBaseMonodCS}{\num{-353780}}

\newcommand{\fmBaseRamServdW}{\num{210736}}

\newcommand{\fmBaseRamSysShare}{73.7}
\newcommand{\fmBaseRamSysPrice}{594}

\newcommand{\fmBaseRamSysdW}{\num{244878}}

\newcommand{\fmBaseSocialPS}{\num{-991605}}
\newcommand{\fmBaseSocialShare}{89.9}
\newcommand{\fmBaseSocialPrice}{56}

\newcommand{\fmBaseSocialdW}{\num{351111}}

\newcommand{\fmHighFuelCompShare}{45.4}
\newcommand{\fmHighFuelCompPrice}{924}
\newcommand{\fmHighFuelCompN}{3}

\newcommand{\fmHighFuelCompdWbase}{\num{-274388}}

\newcommand{\fmHighFuelMonodW}{\num{-8265}}

\newcommand{\fmHighFuelMonodWbase}{\num{-74308}}

\newcommand{\fmHighFuelSocialShare}{89.1}

\newcommand{\fmHighDemandCompShare}{66.0}
\newcommand{\fmHighDemandCompPrice}{838}
\newcommand{\fmHighDemandCompN}{5}

\newcommand{\fmHighDemandMonoPrice}{1025}

\newcommand{\fmHighDemandMonodW}{\num{-456840}}

\newcommand{\fmHighDemandRamSysdW}{\num{240419}}

\newcommand{\fmLowCapacityCompN}{4}

\newcommand{\fmLowCapacityCompLF}{0.953}

\newcommand{\fmLowCapacityMonoPrice}{1026}

\newcommand{\fmLowCapacityMonodWbase}{\num{-6132}}

\newcommand{\fmLowCapacityRamServdWbase}{\num{-9802}}

\newcommand{\fmLowCapacityRamSysdWbase}{\num{-709}}

\newcommand{\fmLowCapacitySocialdWbase}{\num{-22074}}
\newcommand{\fmBaseCompContested}{5}

\newcommand{\fmBaseRamseyShareOfPlannerGain}{70}
\newcommand{\fmBaseRamseyGap}{\num{34142}}
\newcommand{\fmBasePriceGapPct}{12.0}
\newcommand{\fmFuelFixedSetWelfareLoss}{\num{116325}}
\newcommand{\fmHighFuelCompWmin}{\num{2659387}}
\newcommand{\fmHighFuelCompWmax}{\num{2722596}}
\newcommand{\fmValCompDev}{\num{1.4e-07}}

\newcommand{\fmBaseMaxGain}{\num{-273.4}}

\newcommand{\fmBaseCandidates}{758}

\newcommand{\fmHighFuelMaxGain}{\num{-851.3}}

\newcommand{\fmNestBaseCompW}{\num{2725430}}

\newcommand{\fmNestBaseCompN}{3}

\newcommand{\fmNestBaseRamSysdW}{\num{390358}}

\newcommand{\fmNestBaseCompMinusMonoW}{\num{0}}
\newcommand{\fmNestBaseRamseyShareOfPlannerGain}{87}

\newcommand{\fmNestHighDemandCompWmax}{\num{3805196}}

\newcommand{\fmNestBaseDistinctEq}{4}

\newcommand{\fmNestHighFuelDistinctEq}{4}

\newcommand{\fmNestHighDemandDistinctEq}{8}

\newcommand{\fmNestLowCapacityDistinctEq}{4}

\newcommand{\fmSensEightCompN}{3}
\newcommand{\fmSensEightCompContested}{1}

\newcommand{\fmSensNineCompN}{4}

 \graphicspath{{generated/}}
 
\title{A Strategic Freight Market Model for Liner Shipping:\\
Pricing, Entry, and Endogenous Service Choice}
\author{Jonas Westin\thanks{Transportation Research Unit, Centre for Regional Science, Ume\aa{} University. E-mail: \href{mailto:jonas.westin@umu.se}{jonas.westin@umu.se}.}}
\date{November 2026}
 
\begin{document}

\maketitle

\begin{abstract}
National freight transport models usually treat prices and service provision as exogenous. This paper models both as outcomes of strategic interaction among operators of scheduled, loop-based freight services. It compares oligopolistic competition, monopoly, Ramsey pricing with service-level or system-wide cost recovery, and welfare-maximising planning. Shippers across multiple segments choose between direct cyclic liner services and an outside option, using logit demand. Operators select from candidate port rotations, vessel classes and speeds, and set origin–destination prices subject to leg capacities.
Under competition, each operator offers one service. Equilibrium requires that no operator can profitably enter, exit or replace its service after prices re-equilibrate. The solution method computes global pricing best responses and tests all unilateral service deviations, pruning only those excluded by an analytical profit bound. All regimes use the same demand and cost structure.
In a synthetic Baltic–North Sea application, multinomial logit competition sustains more overlapping services and lower prices than monopoly. System-wide Ramsey pricing captures most of the planner's welfare gain without a deficit. Fuel-price and demand shocks affect competitive prices through changes in the service set. When maritime services are close substitutes, entrants avoid overlapping markets and competition can reproduce monopoly outcomes. The results illustrate mechanisms rather than forecast actual markets, showing how market organisation shapes service provision, prices, modal shares and welfare.
\end{abstract}




\newpage
\section{Introduction}
\label{sec:introduction}

Freight transport models usually take transport prices and service supply as given. Operating costs determine which services are feasible, but market prices also reflect demand elasticity, service quality, capacity scarcity and competitive conduct, and the set of services offered is itself the outcome of operators' decisions. This paper develops a strategic freight-market model for loop-based scheduled transport services in which both prices and service provision are endogenous. The work is part of the HERMES project (High Resolution Freight Network Assignment for Swedish Conditions).

In the model, shippers allocate exogenous origin--destination (OD) demand across differentiated liner services and an outside option. Operators set OD-specific prices based on leg capacities and may enter, exit, or replace their service with another candidate service. The same demand, cost, and capacity system is solved under competition, monopoly, Ramsey pricing, and a social planner, making the regimes directly comparable.

The model is applied to a synthetic Baltic--North Sea network. The synthetic setting isolates the economic mechanisms of interest and allows transparent comparisons across regimes and scenarios. The numerical results are model experiments, not estimates of observed freight rates, modal shares or operator behaviour.

\subsection{Background}

Freight transport is frequently organised in loops, rotations, and round trips rather than as independent point-to-point movements. A vehicle or vessel may serve several loading and unloading locations within a single continuous cycle, while different freight flows share capacity along different parts of the same route. Consolidation, backhaul, and repeated service patterns can increase vehicle utilisation and create economies of scale. These structures occur in both local distribution systems and international liner shipping. Empirical work on vessel movements in Swedish waters, for example, documents recurring multi-port patterns in which the same vessel loads and unloads at several calls \parencite{johansson2021}.

Representing such movements explicitly is important because the profitability and capacity use of one freight relation depend on the other flows carried by the same service. Loop-based network representations can therefore describe consolidation, return movements, and shared capacity more naturally than models in which each relation is treated as a separate direct shipment \parencite{asplundwestin2017,westin2024loops,westin2026loops}. At the same time, the number of feasible loops grows rapidly with the number of locations, vessel classes, speeds, and schedules. This creates a direct connection between the economic problem of service provision and the computational literature on liner-shipping network design \parencite{brouer2014,brouer2014matheuristic,koza2020,christiansen2020}.

National freight models provide detailed representations of commodity flows, shipment size, transport chains, terminal use, mode choice, and network assignment. The Swedish Samgods system, for example, follows an aggregate--disaggregate--aggregate structure in which aggregate production--consumption flows are disaggregated into firm-to-firm relations and assigned to logistics alternatives based on transport and inventory costs \parencite{benakivadejong2013,dejongbaak2023,trafikverket2026}. Stochastic extensions improve the behavioural representation by allocating demand probabilistically across logistics alternatives rather than assigning all flow to the minimum-cost option \parencite{abate2019,lindgren2019}.

These advances primarily concern the demand and logistics-choice side of the model. Transport-service supply, freight rates, and strategic interaction among operators are usually represented more simply. This is a consequential limitation when operators serve overlapping markets. Operating costs determine whether a service is technically and commercially feasible, but market prices also reflect demand elasticity, service quality, capacity scarcity, ownership, and competitive conduct. A change in the price, speed, capacity, or route of one operator alters the demand and profitability of competing services. Market outcomes therefore depend not only on demand, geography, and engineering costs, but also on the strategic market decisions through which services and prices are formed.

Scheduled liner shipping provides a useful setting for studying these mechanisms. Liner services combine large fixed operating costs, discrete route choices, heterogeneous transit times, shared leg capacity, and overlapping origin--destination (OD) markets. Advanced network-design models represent cyclic services, vessel deployment, speed, cargo routing, timetables, and transhipment in substantial operational detail \parencite{agarwal2008,brouer2014,meng2014,koza2020}. A separate line of literature examines strategic competition in freight rates, capacity, and service \parencite{wang2014,shah2012,song2023,lyu2023}. Yet national freight-demand models, strategic pricing models, and detailed service-network models have largely developed as separate research traditions. The joint representation of shipper choice, OD-specific strategic pricing, leg capacity, service entry and exit, route replacement, and welfare benchmarks remains comparatively uncommon.

\subsection{Research questions and contributions}

The purpose of the paper is to develop a service-based model of strategic interaction among liner-shipping operators and to assess how operator pricing and service choices affect freight allocation, market structure, profitability, and social welfare.

The paper addresses three main questions:

\begin{enumerate}[leftmargin=*,label=\arabic*.]
  \item How can a loop-based representation of scheduled freight services be combined with discrete-choice demand and strategic operator pricing?
  \item How do endogenous entry, exit, and route replacement affect service provision, prices, freight flows, capacity use, and operator profitability?
  \item How do the resulting market outcomes compare with monopoly, Ramsey pricing under cost recovery, and a social-planner benchmark, and how robust are those outcomes to alternative market conditions, substitution patterns and solution procedures?
\end{enumerate}

The paper makes three contributions.

\begin{enumerate}[leftmargin=*,label=(\roman*)]
  \item It states a two-stage game in which single-service operators choose a loop-based service from a candidate pool and then set OD-specific prices under shared leg capacities, while shippers respond through segmented logit demand with an outside option. A nested-logit variant lets maritime services be closer substitutes for each other than for the outside option (\cref{sec:model}).
  \item It defines a candidate-set equilibrium and a solve--verify--repair procedure that certifies it. Best responses are computed globally, and verification tests every unilateral exit, entry and route replacement in the candidate set against the post-deviation price equilibrium. A deviation is skipped only when a closed-form bound shows that it cannot be profitable, so the certificate does not depend on a screening heuristic (\cref{sec:properties,sec:solution}).
  \item It compares competition with monopoly, Ramsey and planner outcomes on a common demand and cost system. It shows which differences follow from the model structure and which arise through endogenous service choice, and how the advantage of competition over monopoly depends on how closely the services substitute for each other (\cref{sec:properties,sec:study}).
\end{enumerate}

The equilibrium claim is deliberately conditional. It applies to the generated candidate set and the stated numerical tolerances, not to every conceivable route in the unrestricted network-design space.

\subsection{Scope}

The scope is narrower than a complete national freight-market model. Total demand by OD relation and commodity group is exogenous, although its allocation across services is endogenous. Shippers choose among direct scheduled maritime services and an outside option; transfers between services are excluded. In the competitive regime, each operator owns one service, and there is no global constraint on the fleet by vessel class. The outside option is represented through utility rather than through an explicit supply model and resource-cost account. Environmental externalities and the marginal cost of public funds are not included in the welfare measure. The monopoly, Ramsey and planner service sets are obtained by best-improvement search over the full add--remove--swap neighbourhood from several starting sets, not by a proof of global optimality; on reduced instances, full enumeration confirms the search.

The numerical application uses a synthetic test network. This makes it possible to isolate the principal economic mechanisms and validate the algorithms in a controlled setting, and it means that the results are methodological illustrations. Empirical application requires observed freight rates, schedules, vessel and port costs, service capacities, shipment choices and calibrated demand parameters.

\subsection{Outline}

\Cref{sec:literature} positions the framework in four literatures. \Cref{sec:model} presents the model and the market game, and \cref{sec:properties} collects properties that hold by construction, closed-form benchmarks used to validate the solver, and the bounds used to prune the searches. \Cref{sec:solution} describes the solve--verify--repair procedure. \Cref{sec:study} reports the computational study, including a robustness analysis with nested-logit demand. \Cref{sec:discussion} discusses interpretation, limitations and next steps, and \cref{sec:conclusions} concludes.

\newpage

\section{Related literature}
\label{sec:literature}

The model draws on four literatures that have developed largely in parallel: national freight-demand models, discrete-choice models of price competition, models of competition and service design in freight and liner shipping, and models of entry and welfare under fixed costs. Each contains building blocks of the present framework, but they are rarely combined into a single model.

\subsection{Freight-demand models}

Freight forecasting has moved from aggregate commodity-flow models towards behaviourally explicit logistics systems that represent shipment size, transport chains, mode and terminal choice, and network assignment \parencite{chow2010,dejong2016}. \Textcite{dejong2007} model shipment size and transport-chain choice jointly. The Swedish Samgods system implements this aggregate--disaggregate--aggregate structure \parencite{benakivadejong2013,dejongbaak2023,trafikverket2026}, and stochastic versions allocate demand probabilistically across chains rather than assigning all flow to the cheapest one \parencite{abate2019,lindgren2019}. Disaggregate chain-choice and joint mode--shipment-size models show heterogeneous and non-linear responses across commodities \parencite{jensen2019,keya2019}.

These models motivate the segmented logit demand and the outside option used here. On the supply side, however, prices, capacities and service networks enter as given tariffs or cost functions, not as the choices of competing firms. The models explain how shippers choose among a given set of alternatives, not how the prices and services in that set arise.

\subsection{Logit demand and Bertrand pricing}

Random-utility models link service attributes and prices to choice probabilities \parencite{mcfadden1974,train2009} and provide a theory of differentiated-product competition \parencite{anderson1992}. \Textcite{berry1995} combine such demand with supply-side first-order conditions, the structure used here for operators' pricing. Computing Bertrand--Nash prices under logit-type demand is well studied \parencite{morrow2011,gallego2014}, and the optimal prices of a multiproduct firm under multinomial and nested logit have closed forms in terms of the Lambert $W$ function \parencite{lihuh2011}, which the bounds in \cref{sec:bounds} use. Existence and uniqueness of price equilibria require conditions on demand, shares, price bounds or costs \parencite{caplin1991,aksoypierson2013}, and network pricing with logit demand is non-convex even when the network is fixed \parencite{gilbert2015}.

Three limitations carry over to the present model. Multinomial logit imposes independence of irrelevant alternatives within each segment; \cref{sec:nested} relaxes this with a nested logit in which maritime services compete more closely with each other than with the outside option. Single-service operators internalise only their own demand response, unlike multiproduct firms. Moreover, because entry and route replacement change the product set, a unique price equilibrium for a given set of services does not imply a unique market structure.

\subsection{Competition and service design in freight and liner shipping}

\Textcite{shah2012} show that frequency and capacity are strategic dimensions of freight service quality and that competition can produce excess capacity relative to the social optimum. In liner shipping, \textcite{wang2014} analyse logit-demand duopoly competition in freight rates and deployment, including entry deterrence. \Textcite{song2023} study duopoly spot pricing on overlapping OD pairs that share vessel capacity, \textcite{lyu2023} consider capacity allocation and pricing between a carrier and a freight forwarder, and alliance studies examine how cooperation changes prices and welfare \parencite{wang2024alliance,lu2025}. \Textcite{tamannaei2021} model intermodal providers competing with road transport. These models typically fix the set of routes and the number of firms.

The operations research literature represents liner services in far greater detail: cyclic routes, fleet deployment, speed, schedules, cargo routing, and transhipment \parencite{agarwal2008,brouer2014,brouer2014matheuristic,koza2020}, as reviewed by \textcite{meng2014} and \textcite{christiansen2020}. \Textcite{karsten2017} add commodity-specific transit-time limits; their `competitive' networks are attractive in cost and transit time, not outcomes of a game among carriers. Models that join design and pricing study one focal operator facing shipper response and a simplified competitor \parencite{tawfik2019,wang2023integrated}, an entrant facing a given incumbent network \parencite{gelareh2010}, or a fixed number of carriers in a sequential game \parencite{angeloudis2016}. \Textcite{lederer1993} shows that network structure and bundling can destroy the existence or uniqueness of price equilibria, which matters once shippers can combine services.

\subsection{Entry, fixed costs and welfare}

Free entry need not deliver the efficient number of firms when fixed costs are duplicated, and entrants take business from incumbents \parencite{mankiw1986}. Empirical entry models relate market size and firm heterogeneity to market structure \parencite{bresnahan1991,berry1992}, dynamic models add investment and path dependence \parencite{ericson1995,aguirregabiria2012}, and discrete entry games can have several equilibria \parencite{ciliberto2009}. The present model is static and reports multiplicity rather than selecting among equilibria.

For the welfare benchmarks, Ramsey pricing gives the second-best rule when prices must cover costs \parencite{ramsey1927,baumol1970}. \Textcite{horcher2021} review how pricing, capacity, cost recovery and subsidies interact in scheduled public transport, where fixed service costs and shared vehicle capacity raise the same issues as in liner shipping.

\subsection{Position of this paper}

\Cref{tab:literature-positioning} compares representative studies. No single predecessor combines segmented discrete-choice demand, strategic OD pricing, leg capacity, a broad candidate service set with entry, exit and route replacement, and common monopoly, Ramsey and planner counterfactuals.

\begin{table}[htbp]
  \centering
  \caption{Positioning of representative literature.}
  \label{tab:literature-positioning}
  \footnotesize
  \setlength{\tabcolsep}{3.5pt}
  \renewcommand{\arraystretch}{1.15}
  \begin{tabularx}{\textwidth}{>{\raggedright\arraybackslash}p{0.2\textwidth}
                               >{\raggedright\arraybackslash}X
                               >{\raggedright\arraybackslash}p{0.12\textwidth}
                               >{\raggedright\arraybackslash}X
                               >{\raggedright\arraybackslash}X}
    \toprule
    Study & Demand response & Strategic prices & Endogenous services & Market structure \\
    \midrule
    \textcite{dejong2007}; \textcite{abate2019} & Shipment size and chain choice & No & Chains chosen by shippers & No carrier decisions \\
    \textcite{brouer2014}; \textcite{koza2020} & Cargo routing & No & Network, schedule and speed & Single carrier \\
    \textcite{shah2012} & Full price including inventory cost & Yes & Frequency and capacity & Fixed carriers \\
    \textcite{wang2014} & Logit shares & Yes & Frequency and vessel capacity & Duopoly; entry deterrence \\
    \textcite{gelareh2010} & Shares from cost and time & Indirect & Entrant's hub network & Entrant versus incumbent \\
    \textcite{angeloudis2016} & Exogenous OD demand & Yes & Fleet and services & Fixed oligopoly \\
    \textcite{tawfik2019}; \textcite{wang2023integrated} & Random utility, outside option & Leader & Service network design & Focal operator \\
    \textcite{song2023} & Price-sensitive OD demand & Dynamic duopoly & Capacity on given routes & Fixed duopoly \\
    This paper & Segmented logit and nested logit, outside option & OD-specific Bertrand & Candidate cyclic services; entry, exit, replacement & Free single-service entry; monopoly, Ramsey, planner \\
    \bottomrule
  \end{tabularx}
\end{table}

The existing literature leaves several important gaps.
Freight-demand models take supply as given; liner network-design models optimise one carrier or a central system; strategic freight models fix the product set; and entry models abstract from transport networks. The present model sits at an intermediate level of detail. It is more strategically explicit than a network design model and less operationally complete than the leading liner network formulations. Because the service set is discrete, it also needs an explicit equilibrium check: a converged price solution or a local-search endpoint does not establish that no profitable unilateral deviation remains.

\newpage
\section{Model framework}
\label{sec:model}

This section presents the model framework, beginning with the economic environment and the representation of scheduled liner services, including routes, operating costs, transit times and capacity constraints. It then specifies shipper demand, operator pricing and the two-stage market game governing entry, exit and service replacement. Finally, it defines the welfare measure and the alternative market regimes used to compare oligopolistic competition, monopoly, Ramsey pricing and social planning.

\subsection{Economic environment}

The model represents a market for scheduled freight services operated on cyclic transport loops. A service combines a geographical rotation, a vessel class, a sailing speed, a weekly operating pattern, and finite capacity on each leg. Shippers choose among the services that can carry their freight directly between an origin and a destination, together with an outside transport alternative. Operators set prices strategically and may enter the market, leave the market, or replace their current service by another candidate service.

The model therefore contains two interdependent choice layers. At the demand layer, shippers allocate a fixed amount of OD demand across differentiated transport products. At the supply layer, operators determine prices and the set of services offered in the market. Their interaction determines freight volumes, capacity utilisation, operator profits, market structure, and welfare.

Let $K$ denote the set of origin--destination (OD) markets, $G$ the set of commodity groups or demand segments, and $S$ the set of candidate scheduled services. Each service $s\in S$ is associated with a cyclic port rotation and a set of directed legs $L_s$. The vessel class used by service $s$ is denoted $f_s\in F$, and its sailing speed by $v_s$. An active market structure is a subset
\begin{equation}
    A\subseteq S.
\end{equation}
For an OD market $k$, let $S_k(A)\subseteq A$ denote the active services that can carry market $k$ directly and satisfy its transit-time requirement.

Conditional on an active set $A$, operators choose OD-specific prices. Shippers then choose among the available services and the outside option. The resulting demand determines service revenues, leg loads, and profits. In the competitive regime, $A$ is endogenous because operators may enter, exit, or replace their service.

\subsection{Services, schedules, and operating costs}

A candidate service $s$ is defined by an ordered cyclic sequence of distinct ports
\begin{equation}
    R_s=(i_1,i_2,\ldots,i_{n_s},i_1),
\end{equation}
together with vessel class $f_s$ and sailing speed $v_s$. Different vessel classes or speeds on the same geographical rotation are treated as distinct transport products because they imply different costs, capacities, and transit times. The same ports called in the opposite order form a different rotation.

Let $d_{ij}$ denote the maritime distance between ports $i$ and $j$. Sailing time on leg $(i,j)$ is
\begin{equation}
    t_{ij,s}^{\mathrm{sail}}
    =
    \frac{d_{ij}}{v_s}.
    \label{eq:sailing-time}
\end{equation}
If $h_i$ is the port time associated with a call at port $i$, the rotation time of service $s$ is
\begin{equation}
    \tau_s^{\mathrm{rot}}
    =
    \sum_{(i,j)\in L_s} t_{ij,s}^{\mathrm{sail}}
    +
    \sum_{i\in R_s} h_i .
    \label{eq:rotation-time}
\end{equation}
The number of vessels required to maintain one weekly service is
\begin{equation}
    m_s
    =
    \left\lceil
    \frac{\tau_s^{\mathrm{rot}}}{168}
    \right\rceil .
    \label{eq:vessels-required}
\end{equation}
Slack in the rotation, $168m_s-\tau_s^{\mathrm{rot}}$ hours, is treated as buffer time: fuel use below is computed at the nominal speed $v_s$.

Fuel consumption is increasing and convex in speed. A convenient approximation is
\begin{equation}
    FC_f(v)=a_f v^3,
    \label{eq:fuel-function}
\end{equation}
where $a_f$ is vessel-class specific. If $FC_f(v)$ is measured per day, weekly sailing fuel consumption is
\begin{equation}
    Q_s^{\mathrm{fuel}}
    =
    \sum_{(i,j)\in L_s}
    a_{f_s}v_s^3
    \frac{t_{ij,s}^{\mathrm{sail}}}{24}.
    \label{eq:fuel-use}
\end{equation}

The fixed weekly cost of offering service $s$ is
\begin{equation}
    C_s
    =
    m_s C^{\mathrm{charter}}_{f_s}
    +
    q^{\mathrm{fuel}}Q_s^{\mathrm{fuel}}
    +
    \sum_{i\in R_s} C_i^{\mathrm{port}}
    +
    C_s^{\mathrm{other}},
    \label{eq:fixed-service-cost}
\end{equation}
where $q^{\mathrm{fuel}}$ is the fuel price, $C^{\mathrm{charter}}_{f_s}$ is the weekly vessel cost, $C_i^{\mathrm{port}}$ is the cost of a port call, and $C_s^{\mathrm{other}}$ collects other service-level fixed costs (zero in the test instance). Volume-dependent handling and transport costs are represented separately by $c^{\mathrm{var}}_{kgs}$ per forty-foot equivalent unit (FFE) for OD market $k$, commodity group $g$, and service $s$.

\subsection{Service availability and transit time}

A service can carry OD market $k$ directly if both the origin $o_k$ and the destination $d_k$ are called at on its rotation. Because a rotation is a cycle of distinct ports, cargo travels forward along the loop from $o_k$ to $d_k$ on a unique path $P_{ks}$, the sequence of legs from $o_k$ to $d_k$. The transit time is the sailing time on that path plus the port time at the intermediate calls,
\begin{equation}
    T_{ks}
    =
    \sum_{(i,j)\in P_{ks}} t_{ij,s}^{\mathrm{sail}}
    +
    \sum_{i\in I_{ks}} h_i ,
    \label{eq:transit-time}
\end{equation}
where $I_{ks}$ is the set of ports called between $o_k$ and $d_k$. Port time at the origin and destination is not included. The service is available to market $k$ if, in addition, $T_{ks}$ does not exceed the maximum admissible transit time $TT_k$. Define
\begin{equation}
    b_{ks}
    =
    \begin{cases}
        1, & \text{if both $o_k$ and $d_k$ are called at by service $s$},\\
        0, & \text{otherwise}.
    \end{cases}
\end{equation}
The feasible service set for OD market $k$ under market structure $A$ is
\begin{equation}
    S_k(A)
    =
    \left\{
        s\in A:
        b_{ks}=1,\;
        T_{ks}\leq TT_k
    \right\}.
    \label{eq:choice-set}
\end{equation}

The model therefore treats a liner service as a direct scheduled transport product: a shipment remains on the same service between origin and destination.

\subsection{Shipper demand}
\label{sec:demand}

Weekly demand $D_{kg}$ is exogenous for each OD market $k$ and commodity group $g$. The model determines how this demand is distributed across the available scheduled services and an outside transport alternative. The outside option represents transport arrangements not explicitly modelled as services in $A$, such as road, rail, another maritime product, or another logistics arrangement.

Let $p_{ks}$ denote the price per FFE charged by service $s$ in OD market $k$. A common price is charged to all commodity groups in the same OD--service pair, while the groups may differ in their sensitivity to price and time.

The systematic utility of service $s$ for group $g$ in market $k$ is
\begin{equation}
    V_{kgs}
    =
    \gamma_{0,g}
    -
    \gamma_{T,g}T_{ks}
    -
    \beta_g p_{ks},
    \qquad
    s\in S_k(A),
    \label{eq:utility}
\end{equation}
where $\beta_g>0$ is the marginal disutility of price and $\gamma_{T,g}>0$ is the marginal disutility of transit time. The systematic utility of the outside option is $V^{\mathrm{out}}_{kg}$.

Under multinomial logit demand, the market share of service $s$ is
\begin{equation}
    \phi_{kgs}(p;A)
    =
    \frac{\exp(V_{kgs})}
    {
        \exp(V^{\mathrm{out}}_{kg})
        +
        \sum_{r\in S_k(A)}
        \exp(V_{kgr})
    },
    \label{eq:logit-share}
\end{equation}
and the outside-option share is
\begin{equation}
    \phi^{\mathrm{out}}_{kg}(p;A)
    =
    \frac{\exp(V^{\mathrm{out}}_{kg})}
    {
        \exp(V^{\mathrm{out}}_{kg})
        +
        \sum_{r\in S_k(A)}
        \exp(V_{kgr})
    }.
\end{equation}
The volume carried by service $s$ is therefore
\begin{equation}
    x_{kgs}(p;A)
    =
    D_{kg}\phi_{kgs}(p;A).
    \label{eq:demand-allocation}
\end{equation}

The logit formulation provides a direct analytical link between operator prices and freight allocation. For two services $s$ and $r$ serving the same OD market,
\begin{equation}
    \frac{\partial \phi_{kgs}}{\partial p_{kr}}
    =
    \begin{cases}
    -\beta_g \phi_{kgs}(1-\phi_{kgs}),
        & r=s,\\[0.4em]
    \beta_g \phi_{kgs}\phi_{kgr},
        & r\neq s.
    \end{cases}
    \label{eq:logit-derivatives}
\end{equation}
Thus, an increase in a service's own price reduces its share, while an increase in a competing service's price raises it. Multinomial logit is the main specification of the paper.

\paragraph{Nested-logit variant.}
Under multinomial logit, a new service draws demand from rival services and from the outside option in proportion to their shares. When maritime services are closer substitutes for each other than for road or rail, this overstates how much a new service expands the maritime market and understates how much it takes from incumbents. The variant used in \cref{sec:nested} places the outside option in a nest of its own and all maritime services of market $k$ in one nest with parameter $\rho_g\in(0,1]$. With
\begin{equation}
    w_{kgs}=\exp\!\left(\frac{V_{kgs}}{\rho_g}\right),
    \qquad
    N_{kg}=\sum_{r\in S_k(A)}w_{kgr},
    \qquad
    \sigma_{kgs}=\frac{w_{kgs}}{N_{kg}},
    \label{eq:nest-sum}
\end{equation}
the market share of service $s$ is the share within the maritime nest times the probability of choosing the nest,
\begin{equation}
    \phi_{kgs}
    =
    \sigma_{kgs}\,
    \frac{N_{kg}^{\rho_g}}{\exp(V^{\mathrm{out}}_{kg})+N_{kg}^{\rho_g}} .
    \label{eq:nested-share}
\end{equation}
The derivatives that replace \eqref{eq:logit-derivatives} are
\begin{equation}
    \frac{\partial \phi_{kgs}}{\partial p_{kr}}
    =
    \begin{cases}
    -\beta_g \phi_{kgs}\left[\dfrac{1}{\rho_g}-\dfrac{1-\rho_g}{\rho_g}\sigma_{kgs}-\phi_{kgs}\right],
        & r=s,\\[0.9em]
    \beta_g \phi_{kgs}\left[\dfrac{1-\rho_g}{\rho_g}\sigma_{kgr}+\phi_{kgr}\right],
        & r\neq s.
    \end{cases}
    \label{eq:nested-derivatives}
\end{equation}
With $\rho_g=1$ the nest sum is $N_{kg}=\sum_r\exp(V_{kgr})$, \eqref{eq:nested-share} reduces to \eqref{eq:logit-share}, and \eqref{eq:nested-derivatives} reduces to \eqref{eq:logit-derivatives}. A smaller $\rho_g$ makes the maritime services closer substitutes: a price cut by one service then draws relatively more of its extra traffic from rival services and less from the outside option. When only one service serves market $k$, $N_{kg}^{\rho_g}=\exp(V_{kgs})$ and the nested and multinomial shares coincide for every $\rho_g$; \cref{prop:uncontested} uses this.

\subsection{Leg-specific capacity}

Freight assigned to a service uses capacity on all legs between its origin and destination. Let
\begin{equation}
    a_{ks\ell}
    =
    \begin{cases}
        1, & \text{if OD market $k$ uses leg $\ell$ of service $s$},\\
        0, & \text{otherwise}.
    \end{cases}
\end{equation}
The load on leg $\ell\in L_s$ is
\begin{equation}
    q_{s\ell}(p;A)
    =
    \sum_{k\in K}
    \sum_{g\in G}
    a_{ks\ell}x_{kgs}(p;A).
    \label{eq:leg-load}
\end{equation}
Capacity feasibility requires
\begin{equation}
    q_{s\ell}(p;A)
    \leq
    u_{s\ell},
    \qquad
    \forall s\in A,\;\ell\in L_s,
    \label{eq:capacity}
\end{equation}
where $u_{s\ell}$ is the weekly capacity of the vessel class of $s$, in FFE.

The capacity constraints couple otherwise separate OD pricing decisions. When capacity is slack, each OD price mainly trades off margin against demand in that market. When one or more legs are binding, an additional unit of demand in one OD market displaces potential freight in all other OD markets using the same constrained leg.

Let $\lambda_{s\ell}\geq0$ denote the shadow value of capacity on leg $\ell$. A useful economic marginal-cost measure for an OD--commodity flow is
\begin{equation}
    MC^{\mathrm{econ}}_{kgs}
    =
    c^{\mathrm{var}}_{kgs}
    +
    \Lambda_{ks},
    \qquad
    \Lambda_{ks}
    =
    \sum_{\ell\in L_s}
    a_{ks\ell}\lambda_{s\ell}.
    \label{eq:economic-mc}
\end{equation}
The shadow-price component $\Lambda_{ks}$ is an opportunity cost rather than a cash expenditure.

\subsection{Operator profit and price competition}

In the competitive regime, each active service is operated by one independent firm. Weekly profit is
\begin{equation}
    \pi_s(p;A)
    =
    \sum_{k\in K}
    \sum_{g\in G}
    x_{kgs}(p;A)
    \left(
        p_{ks}-c^{\mathrm{var}}_{kgs}
    \right)
    -
    C_s.
    \label{eq:profit}
\end{equation}
For a fixed active set $A$, operator $s$ chooses
\begin{equation}
    p_s
    =
    \left(p_{ks}\right)_{k:s\in S_k(A)}
\end{equation}
subject to common lower and upper price bounds and the leg-capacity constraints.

Because demand depends on competitors' prices, the feasible capacity set of operator $s$ also depends on competitors' decisions. The capacitated price game can therefore be interpreted as a generalised Nash game \parencite{facchinei2007}. If the capacity constraints are slack, it reduces to a conventional differentiated-product Bertrand game.

For fixed competitor prices $p_{-s}$, the best-response correspondence is
\begin{equation}
    BR_s(p_{-s};A)
    =
    \arg\max_{p_s\in\mathcal P_s(p_{-s};A)}
    \pi_s(p_s,p_{-s};A),
    \label{eq:best-response}
\end{equation}
where $\mathcal P_s(p_{-s};A)$ is the set of prices satisfying the price bounds and capacity constraints.

A price equilibrium for active set $A$ is a vector $p^*(A)$ such that
\begin{equation}
    p_s^*(A)
    \in
    BR_s(p_{-s}^*(A);A)
    \qquad
    \forall s\in A.
    \label{eq:price-equilibrium}
\end{equation}
The corresponding equilibrium profit is
\begin{equation}
    \pi_s^*(A)
    =
    \pi_s(p^*(A);A).
\end{equation}

For an unconstrained OD market, the own-price derivative of profit is
\begin{equation}
    \frac{\partial\pi_s}{\partial p_{ks}}
    =
    \sum_{g\in G}
    \left[
        x_{kgs}
        +
        \left(p_{ks}-c^{\mathrm{var}}_{kgs}\right)
        \frac{\partial x_{kgs}}{\partial p_{ks}}
    \right].
    \label{eq:profit-gradient}
\end{equation}
With binding leg capacity, $c^{\mathrm{var}}_{kgs}$ is replaced by the economic marginal cost \eqref{eq:economic-mc} in the first-order conditions.

\subsection{The market game: entry, exit, and route replacement}
\label{sec:game}

Price equilibrium is conditional on the active set. The active set is itself determined in a two-stage game.

\paragraph{Players and strategies.}
There is a set $N$ of potential operators, large enough that every candidate service could be operated. In the first stage, each operator $i\in N$ chooses $a_i\in S\cup\{\emptyset\}$, where $\emptyset$ denotes staying out. Each candidate service can be operated by at most one operator, so the active set $A=\{a_i:a_i\neq\emptyset\}$ contains no duplicates. This rules out \emph{clone entry}, in which a second operator runs a copy of an active service; \cref{sec:multiplicity} discusses this choice.

\paragraph{Prices and payoffs.}
In the second stage, the operators of the active services observe $A$ and set prices simultaneously, as in \eqref{eq:best-response}--\eqref{eq:price-equilibrium}. If the price stage has several equilibria for some $A$, $p^*(A)$ denotes the one selected by the solution procedure, which starts from the stand-alone prices defined in \cref{sec:price-stage}. An active operator earns $\pi_s^*(A)$; an operator that stays out earns zero.

\paragraph{Unilateral deviations.}
Given the price-stage outcome $p^*(\cdot)$, an operator can deviate in the first stage in three ways, and every deviation is evaluated after all active operators have re-equilibrated prices. An active operator may exit. Exit yields zero, so an active service with $\pi_s^*(A)<0$ has a profitable exit deviation. An inactive operator may enter with an unoccupied candidate $r\in S\setminus A$. Define
\begin{equation}
    A^{+r}
    =
    A\cup\{r\}.
\end{equation}
Entry is profitable if the entrant's profit is positive:
\begin{equation}
    \Delta\pi_r^{\mathrm{entry}}
    =
    \pi_r^*(A^{+r})
    >
    0.
    \label{eq:entry-gain}
\end{equation}
An active operator may instead replace its service $s\in A$ by an unoccupied candidate $r\notin A$. Let
\begin{equation}
    A^{s\rightarrow r}
    =
    \left(A\setminus\{s\}\right)\cup\{r\}.
\end{equation}
The gain from replacement is
\begin{equation}
    \Delta\pi_{s\rightarrow r}^{\mathrm{switch}}
    =
    \pi_r^*(A^{s\rightarrow r})
    -
    \pi_s^*(A).
    \label{eq:switch-gain}
\end{equation}

\paragraph{Equilibrium.}
A candidate-set competitive equilibrium is an active set $A^*$ and associated price equilibrium $p^*(A^*)$ satisfying
\begin{align}
    \pi_s^*(A^*) &\geq -\varepsilon_\pi
    && \forall s\in A^*,
    \label{eq:eq-exit}\\
    \pi_r^*(A^*\cup\{r\}) &\leq \varepsilon_\pi
    && \forall r\in S\setminus A^*,
    \label{eq:eq-entry}\\
    \pi_r^*(A^{*\,s\rightarrow r})-\pi_s^*(A^*)
    &\leq \varepsilon_\pi
    && \forall s\in A^*,\;r\in S\setminus A^*,
    \label{eq:eq-switch}
\end{align}
where $\varepsilon_\pi\geq0$ is a numerical tolerance.
With $\varepsilon_\pi=0$, conditions \eqref{eq:eq-exit}--\eqref{eq:eq-switch} are the first-stage Nash conditions given the selected price-stage play. If a price equilibrium exists for every active set, so that $p^*(\cdot)$ is defined in every subgame, $(A^*,p^*)$ is a pure-strategy subgame-perfect equilibrium of the two-stage game restricted to the candidate set $S$. Deviations excluded by the bound \eqref{eq:standalone-bound} are unprofitable under any rival prices, so this conclusion does not depend on which price equilibrium is selected in the subgames they lead to.

\subsection{Welfare}

Consumer surplus for OD market $k$ and group $g$ follows from the logsum:
\begin{equation}
    CS_{kg}(p;A)
    =
    \frac{D_{kg}}{\beta_g}
    \ln
    \left[
        \exp(V^{\mathrm{out}}_{kg})
        +
        \sum_{s\in S_k(A)}
        \exp(V_{kgs})
    \right].
    \label{eq:consumer-surplus-market}
\end{equation}
Under the nested-logit variant, the sum in brackets is replaced by $N_{kg}^{\rho_g}$ from \eqref{eq:nest-sum}:
\begin{equation}
    CS_{kg}(p;A)
    =
    \frac{D_{kg}}{\beta_g}
    \ln
    \left[
        \exp(V^{\mathrm{out}}_{kg})
        +
        N_{kg}^{\rho_g}
    \right].
    \label{eq:consumer-surplus-nested}
\end{equation}
Aggregate consumer surplus is
\begin{equation}
    CS(p;A)
    =
    \sum_{k\in K}\sum_{g\in G} CS_{kg}(p;A).
    \label{eq:consumer-surplus}
\end{equation}
Producer surplus before any subsidy is
\begin{equation}
    PS(p;A)
    =
    \sum_{s\in A}\pi_s(p;A),
    \label{eq:producer-surplus}
\end{equation}
and total welfare is
\begin{equation}
    W(p;A)
    =
    CS(p;A)+PS(p;A).
    \label{eq:welfare}
\end{equation}

The outside alternative enters through utility and therefore affects shipper choice and consumer surplus, but its resource costs are not separately represented. The planner comparisons should consequently be interpreted as welfare comparisons within the maintained freight-market system.

Consumer surplus in \eqref{eq:consumer-surplus-market} depends on how utilities are normalised. Adding a constant $c$ to every utility in market $k$ and segment $g$, including the outside option, leaves all shares unchanged and shifts $CS_{kg}$ by $D_{kg}c/\beta_g$. Welfare differences between regimes with the same demand are invariant to the normalisation; welfare levels and percentage changes in them are not (\cref{prop:normalisation}). The results in \cref{sec:study} are therefore reported as differences, and the normalisation used is stated in \cref{app:parameters}.

\subsection{Alternative market regimes}

The same demand system, service technology, costs, and capacities are used to compare four regimes.

\subsubsection{Competition}

Competition is defined by the price and service-structure equilibrium in \eqref{eq:price-equilibrium}--\eqref{eq:eq-switch}. Each active service is assumed to be operated by an independent firm.

\subsubsection{Monopoly}

Under monopoly, one owner controls all active services and internalises substitution among them. Total profit is
\begin{equation}
    \Pi^{M}(p;A)
    =
    \sum_{s\in A}\pi_s(p;A).
\end{equation}
The monopoly problem is
\begin{equation}
    \max_{A\subseteq S,\;p}
    \;
    \Pi^{M}(p;A)
    \label{eq:monopoly}
\end{equation}
subject to price bounds and all leg-capacity constraints.

\subsubsection{Ramsey pricing}

We consider two versions of Ramsey pricing that maximise social welfare subject to different cost-recovery requirements. Both determine prices and the set of active services endogenously, but differ in whether financial surpluses can be redistributed across services.

Under system-wide cost recovery, the service network must break even in aggregate, allowing profitable services to cross-subsidise unprofitable ones. The Ramsey problem is
\begin{align}
    \max_{A\subseteq S,\;p}
    \quad & W(p;A)
    \label{eq:ramsey-system-problem}\\
    \text{s.t.}\quad
    & \sum_{s\in A}\pi_s(p;A)\geq0,
    \nonumber
\end{align}
together with price and capacity constraints.

Under service-level cost recovery, each active service must cover its own costs, ruling out cross-subsidisation between services. The corresponding problem is
\begin{align}
    \max_{A\subseteq S,\;p}
    \quad & W(p;A)
    \label{eq:ramsey-service-problem}\\
    \text{s.t.}\quad
    & \pi_s(p;A)\geq0
    \qquad
    \forall s\in A.
    \nonumber
\end{align}
Comparing these two regimes isolates the welfare implications of allowing cross-subsidisation across services while maintaining overall financial viability.

\subsubsection{Social planner}

The social-planner benchmark chooses the active services and prices to maximise welfare without a financial cost-recovery constraint:
\begin{equation}
    \max_{A\subseteq S,\;p}
    W(p;A)
    \label{eq:social-planner}
\end{equation}
subject to price bounds and leg capacities. When the price bounds do not bind and variable costs do not differ across segments, the planner sets each price equal to economic marginal cost \eqref{eq:economic-mc}; otherwise prices satisfy a weighted version of this condition (\cref{sec:benchmarks}). The planner's deficit therefore equals the fixed costs of the active services net of any capacity rents. The planner finances fixed costs through transfers; it does not need prices below marginal cost.

\newpage
\section{Structural properties and analytical bounds}
\label{sec:properties}

This section collects results that follow from the structure of the model rather than from the computations. They serve three purposes. The properties in \cref{sec:structural} separate structural results from findings that depend on the endogenous service set, and they provide checks that the numerical procedures must pass. The closed-form prices in \cref{sec:benchmarks} are benchmarks for validating the solvers. The bounds in \cref{sec:bounds} let the searches in \cref{sec:solution} skip deviations and service sets that cannot be profitable. Throughout, an optimum of a regime means a global optimum of its problem over $A\subseteq S$ and the admissible prices. All results hold for multinomial and nested logit alike.


\subsection{Structural properties}
\label{sec:structural}

\begin{property}[Welfare ordering]
\label{prop:ordering}
Let $W^{C}$ be welfare at a candidate-set competitive equilibrium with $\varepsilon_\pi=0$, and let $W^{M}$, $W^{R,\mathrm{serv}}$, $W^{R,\mathrm{sys}}$ and $W^{\SO}$ be welfare at the optima of \eqref{eq:monopoly}, \eqref{eq:ramsey-service-problem}, \eqref{eq:ramsey-system-problem} and \eqref{eq:social-planner}. Then
\begin{equation}
    W^{\SO}
    \;\geq\;
    W^{R,\mathrm{sys}}
    \;\geq\;
    \max\left\{W^{R,\mathrm{serv}},\,W^{M}\right\},
    \qquad
    W^{R,\mathrm{serv}}
    \;\geq\;
    W^{C}.
    \label{eq:welfare-ordering-property}
\end{equation}
\end{property}

A competitive equilibrium satisfies $\pi_s\geq0$ for every active service by \eqref{eq:eq-exit}, so it is feasible for the service-level Ramsey problem. The monopoly optimum earns $\Pi^M\geq0$, because the empty service set earns zero, so it is feasible for the system-wide Ramsey problem. The service-level constraint set is contained in the system-wide one, which is contained in the planner's, and all regimes share the same price bounds and capacities.
The ordering therefore holds by construction. Because the service-set searches
are started from the solutions they must dominate (\cref{sec:planner-search}), a
reported violation would show that a fixed-set problem was not solved correctly,
which makes \eqref{eq:welfare-ordering-property} a check on the solvers rather
than a finding.
The property is silent on the comparison between competition and monopoly, which is the informative one.

\begin{property}[Fuel price]
\label{prop:fuel}
The fuel price $q^{\mathrm{fuel}}$ enters only the fixed service cost \eqref{eq:fixed-service-cost}. For a given active set $A$, the competitive price equilibria and the monopoly and planner optimal prices are therefore independent of $q^{\mathrm{fuel}}$. Ramsey prices depend on $q^{\mathrm{fuel}}$ whenever the cost-recovery constraint binds.
\end{property}

The fixed cost $C_s$ enters profit \eqref{eq:profit} and welfare \eqref{eq:welfare} additively, so it does not affect which prices maximise them; it does enter the Ramsey constraints. Because speed is attached to each candidate service, a fuel-price change can still alter competitive and monopoly prices, but only by changing the active set, for example towards a slower variant of the same rotation.

\begin{property}[Demand scale]
\label{prop:demand}
Suppose every $D_{kg}$ is multiplied by the same factor $\lambda>0$ and no leg capacity binds before or after the change. For a given active set, the competitive price equilibria and the monopoly and planner optimal prices are unchanged. Consumer surplus and variable profit scale by $\lambda$; fixed costs do not.
\end{property}

With slack capacity, operator profit is $\lambda$ times variable profit minus $C_s$, the monopoly objective is $\lambda$ times total variable profit minus total fixed cost, and the planner objective is $\lambda$ times the sum of consumer surplus and variable profit minus total fixed cost. None of these maximisers depends on $\lambda$. Under uniform demand growth, a change in competitive or monopoly prices for an unchanged set of services must therefore come from binding capacity. Otherwise prices change only through entry, exit or route replacement. Welfare also scales mechanically: at given prices and services it rises by $(\lambda-1)$ times consumer surplus plus variable profit. The Ramsey constraints change with $\lambda$, so Ramsey prices do respond.

\begin{property}[Utility normalisation]
\label{prop:normalisation}
Adding a constant to all utilities in a market--segment pair leaves all shares, prices, profits and service choices unchanged and shifts $CS_{kg}$ by $D_{kg}c/\beta_g$. Welfare differences between regimes or scenarios with the same demand are invariant to the normalisation. Welfare levels, percentage changes, and comparisons across scenarios with different demand are not.
\end{property}

The shift does not depend on prices or on the active set, so it cancels from every comparison that holds demand fixed, including the Ramsey and planner choices among service sets. Under nested logit the same holds, because adding $c$ to every utility multiplies $N_{kg}^{\rho_g}$ and $\exp(V^{\mathrm{out}}_{kg})$ by the same factor $e^{c}$.

\begin{property}[Uncontested markets]
\label{prop:uncontested}
Let $A$ be an active set in which every OD market is served by at most one active service. Then (a) shares, prices, profits and welfare for $A$ are the same under nested logit, for every $\rho_g\in(0,1]$, as under multinomial logit; and (b) the competitive price equilibrium for $A$ maximises total profit over the prices for $A$, so competitive and monopoly prices, profits and welfare for $A$ coincide.
\end{property}

Part (a) follows from \eqref{eq:nested-share}: with one service in the nest, $\sigma_{kgs}=1$ and $N_{kg}^{\rho_g}=\exp(V_{kgs})$. For part (b), the profit of each service then depends only on its own prices, because no rival serves its markets and its capacity constraints involve only its own legs. Total profit is a sum of functions of disjoint price vectors, and maximising the sum is the same as maximising each term, which is what a best response does. The property concerns a given set. Whether competition and monopoly select the same set is an empirical question for the instance, which \cref{sec:nested} answers.

\subsection{Closed-form prices}
\label{sec:benchmarks}

With a single demand segment, slack capacity and interior prices, profit separates across OD markets, because the shares in market $k$ depend only on prices in $k$. The pricing conditions then have closed forms that the numerical solvers must reproduce. Write $m_{ks}=p_{ks}-c^{\mathrm{var}}_{ks}$ for the mark-up of service $s$ in market $k$ and $\phi_{k0}$ for the outside share. For a single-service operator, the first-order condition \eqref{eq:profit-gradient} with the derivatives \eqref{eq:nested-derivatives} gives
\begin{equation}
    m_{ks}
    =
    \frac{1}{\beta\left[\dfrac{1}{\rho}-\dfrac{1-\rho}{\rho}\,\sigma_{ks}-\phi_{ks}\right]},
    \label{eq:markup-competition}
\end{equation}
which is $1/[\beta(1-\phi_{ks})]$ under multinomial logit ($\rho=1$). A smaller $\rho$ lowers the competitive mark-up of services that share a market with rivals, because their demand is more price-elastic, and leaves the mark-up of a service that is alone in its market ($\sigma_{ks}=1$) unchanged. A monopolist that owns all services in market $k$ sets the same mark-up on each of them,
\begin{equation}
    m_{ks}
    =
    \frac{1}{\beta\,\phi_{k0}},
    \label{eq:markup-monopoly}
\end{equation}
for every $\rho$, because with a common mark-up the nest behaves as one product facing the outside option \parencite{lihuh2011}. The planner sets $m_{ks}=0$. With binding legs, $c^{\mathrm{var}}$ is replaced by the economic marginal cost \eqref{eq:economic-mc}. \Cref{sec:validation} reports the largest deviation of the numerical solutions from these values.

\subsection{Profit and welfare bounds}
\label{sec:bounds}

The searches in \cref{sec:solution} need upper bounds on what a set of services can earn or achieve, so that deviations and neighbouring sets that cannot improve on the current best can be skipped without solving them. Two relaxations give closed forms: drop the leg capacities, and let each segment pay its own price. Under these relaxations, in each market--segment row $(k,g)$ an optimal price vector has a common mark-up over the variable cost, and the row is described by one number,
\begin{equation}
    a_{kg}(A)
    =
    \exp\!\left(-V^{\mathrm{out}}_{kg}\right)
    \left[
    \sum_{s\in S_k(A)}
    \exp\!\left(\frac{\gamma_{0,g}-\gamma_{T,g}T_{ks}-\beta_g c^{\mathrm{var}}_{kgs}}{\rho_g}\right)
    \right]^{\rho_g},
    \label{eq:row-aggregate}
\end{equation}
the attractiveness of the maritime nest at marginal-cost prices relative to the outside option. With the normalised mark-up $y=\beta_g m$, the outside share is $1/(1+a_{kg}e^{-y})$, consumer surplus is $(D_{kg}/\beta_g)[V^{\mathrm{out}}_{kg}+\ln(1+a_{kg}e^{-y})]$, and variable profit is $(D_{kg}/\beta_g)\,y\,(1-\phi_{kg0})$. Maximising each objective row by row gives, with $W(\cdot)$ the principal branch of the Lambert $W$ function and $F(A)=\sum_{s\in A}C_s$,
\begin{align}
    \overline{W}^{\SO}(A)
    &=
    \sum_{k,g}\frac{D_{kg}}{\beta_g}\left[V^{\mathrm{out}}_{kg}+\ln\bigl(1+a_{kg}(A)\bigr)\right]-F(A),
    \label{eq:bound-social}\\
    \overline{\Pi}^{M}(A)
    &=
    \sum_{k,g}\frac{D_{kg}}{\beta_g}\,W\!\left(\frac{a_{kg}(A)}{e}\right)-F(A).
    \label{eq:bound-monopoly}
\end{align}
The planner's bound \eqref{eq:bound-social} prices at marginal cost ($y=0$), and the monopoly bound \eqref{eq:bound-monopoly} uses the optimal normalised mark-up $y=1+W(a/e)$ \parencite{lihuh2011}. For system-wide Ramsey pricing, the mark-up in each row is $y(\theta)=\theta+W(\theta a_{kg}e^{-\theta})$, where $\theta=\mu/(1+\mu)\in[0,1)$ is the normalised budget multiplier and is chosen so that total variable profit equals $F(A)$; the bound $\overline{W}^{R}(A)$ is the resulting consumer surplus, or $-\infty$ if $\overline{\Pi}^{M}(A)<0$. Because the service-level constraints are tighter, $\overline{W}^{R}(A)$ also bounds service-level Ramsey welfare.

Each bound is an upper bound on the optimum of the corresponding problem for the set $A$, and it is exact for one segment with slack capacity. For a single service $r$, $B_r=\overline{\Pi}^{M}(\{r\})$ also bounds the profit that $r$ can earn in any active set,
\begin{equation}
    \pi_r^*(A\cup\{r\})\;\leq\;B_r
    \qquad
    \text{for every }A\subseteq S,
    \label{eq:standalone-bound}
\end{equation}
because rival services can only lower the demand for $r$ at every price: the nest sum $N_{kg}$ rises with every rival, and $\phi_{kgr}$ in \eqref{eq:nested-share} decreases in $N_{kg}$ for given $w_{kgr}$. A candidate with $B_r\leq0$ cannot break even in any market structure.

\subsection{Existence, uniqueness, and multiplicity}
\label{sec:multiplicity}

For single-service operators with one logit segment and slack capacity, profit separates across OD markets, and a price equilibrium exists and is unique \parencite{caplin1991,anderson1992}. Two features of the model can break this. First, when segments with different price sensitivities face a common price, demand for a service is a mixture of logits and profit need not be quasi-concave in the service's own price: an operator may prefer to serve both segments at a low price or only the time-sensitive segment at a high price. \Textcite{aksoypierson2013} give conditions for existence and uniqueness under mixed logit demand. Second, binding leg capacities make each operator's feasible prices depend on rivals' prices, which turns the price stage into a generalised Nash equilibrium problem that may have many equilibria \parencite{facchinei2007}. The solution procedure therefore computes best responses globally and uses several starting prices to detect multiple price equilibria (\cref{sec:solution}).

The service stage raises separate issues. Discrete entry games with fixed costs can have several equilibria \parencite{ciliberto2009}, and because the first-stage game is not in general a potential game, best-improvement dynamics over entry and replacement moves can cycle. Finally, the exclusion of clone entry in \cref{sec:game} is a modelling choice. Under multinomial logit, a copy of an active service attracts positive demand, so allowing clones would add one entry deviation per active service. \Cref{app:diagnostics} reports whether any such clone would be profitable in the reported equilibria.

Existence of a pure-strategy equilibrium is therefore not guaranteed. At the price stage it can fail when several segments face a common price or when capacities bind, and at the service stage a cyclic pattern of profitable entries and replacements can rule out any stable service set; the finite first stage then has only mixed-strategy equilibria. The solution procedure does not assume existence: it reports a cycle when the search returns to a set it has visited and counts every price stage that fails to converge. Neither occurred in the computations of \cref{sec:study}.

\newpage
\section{Solution method}
\label{sec:solution}

This section describes how the equilibria and optima of \cref{sec:model} are computed. The problem is split into a pricing stage for a fixed set of services and an outer search over service sets. Every candidate equilibrium is then verified against all unilateral deviations, using the bounds of \cref{sec:bounds} to skip deviations that cannot be profitable.

\subsection{Decomposition principle}

The model combines a continuous strategic pricing problem with discrete service-structure decisions. Solving the entire problem in one monolithic formulation would require simultaneous treatment of logit demand, capacity constraints, fixed service costs, discrete route choices, and non-cooperative operator objectives. The solution method therefore exploits the natural decomposition between a \emph{price subproblem} and a \emph{service-structure problem}.

For a fixed active set $A$, the inner problem is to determine the price equilibrium $p^*(A)$. The outer problem modifies $A$ through entry, exit, and route replacement. Every discrete deviation is evaluated using the post-deviation price equilibrium rather than holding incumbent prices fixed. This is essential because the profitability of a new or modified service depends on the strategic price response of all active operators.

The same decomposition carries over to the alternative market regimes. Competition uses a price-equilibrium subproblem, whereas monopoly, Ramsey, and the social planner solve a continuous optimisation problem conditional on $A$. In every regime, the discrete service-set problem is handled by an outer search, and the closed-form bounds of \cref{sec:bounds} decide which discrete moves need to be solved at all.

\subsection{Price equilibrium for a fixed service set}
\label{sec:price-stage}

For a fixed active set $A$, the price equilibrium is computed by sequential best responses. Let $p^{(n)}$ denote the complete price vector at iteration $n$. In a Gauss--Seidel sweep, operator $s$ responds to the most recent prices of operators already updated in the current sweep and the previous prices of the remaining operators, with damping
\begin{equation}
    p_s^{(n+1)}
    =
    (1-\theta)p_s^{(n)}
    +
    \theta
    BR_s\!\left(p_{-s}^{(n+1,n)};A\right),
    \qquad
    0<\theta\leq1,
    \label{eq:damped-br}
\end{equation}
where $p_{-s}^{(n+1,n)}$ denotes the mixed vector containing the most recent available prices of competitors. The iteration uses $\theta=1$ and is repeated with $\theta=0.5$ if it does not converge. It starts from the \emph{stand-alone prices}: each service's optimal price in each market if it faced only the outside option, ignoring capacity. This start defines the selected equilibrium $p^*(A)$ when several exist (\cref{sec:game}).

\paragraph{Global best responses.}
For fixed rival prices, the shares in market $k$ depend only on prices in $k$, so the profit of operator $s$ is a sum of one-dimensional functions, one per market it serves; only the leg capacities couple them. For multipliers $\lambda_s\geq0$ on the legs of $s$, the Lagrangian
\begin{equation}
    \mathcal L_s(p_s,\lambda_s)
    =
    \sum_{k}\sum_{g}x_{kgs}(p_{ks})\left(p_{ks}-c^{\mathrm{var}}_{kgs}-\Lambda_{ks}\right)
    +\sum_{\ell\in L_s}\lambda_{s\ell}u_{s\ell}-C_s
    \label{eq:br-lagrangian}
\end{equation}
separates into one-dimensional problems in which each market is priced against the economic marginal cost \eqref{eq:economic-mc}. Each one-dimensional problem is solved globally: the objective is evaluated on a grid of 160 prices that is geometric in the mark-up over the lowest variable cost, and safeguarded Newton steps refine the best grid cell. If the resulting prices satisfy all leg capacities, they are the best response. Otherwise the multipliers of the binding legs are found by a projected Newton method on the capacity conditions, with leg-by-leg bisection as a fallback.

By the Lagrangian sufficiency theorem \parencite{everett1963}, prices that maximise \eqref{eq:br-lagrangian} for some $\lambda_s\geq0$, satisfy the capacity constraints, and leave no slack on legs with $\lambda_{s\ell}>0$ are a global best response. With one segment, revenue is concave in the quantity sold under logit, so such multipliers exist. With two segments and a binding leg, the computed response is optimal for the computed multipliers and feasible, but it may fall short of the best response.
%
Under multinomial logit, no leg binds in any of the competitive equilibria reported in \cref{sec:study}. Under nested logit with halved capacity, one leg binds (\cref{tab:nest-diagnostics}); that equilibrium has no contested market, so by \cref{prop:uncontested} each best response is independent of rival prices, and the monopoly solver of
\cref{sec:planner-search} returns the same prices.
This agreement between two independent methods supports, but does not guarantee, global optimality of the computed response.

A global search matters because profit under logit demand is not concave in price. It is nearly flat at high prices, where demand and its derivative are close to zero, so a local method started there can stop at a price that serves almost no one. An `equilibrium' found this way can contain services that earn almost nothing, which in turn can make entry look unprofitable. Evaluating every best response on a grid, refined by Newton steps and checked on a finer grid, removes both failure modes up to the grid resolution.

\paragraph{Convergence and certificate.}
Two convergence measures are used. The first is the maximum price change,
\begin{equation}
    r_p^{(n)}
    =
    \max_{s\in A}
    \left\|
        p_s^{(n+1)}-p_s^{(n)}
    \right\|_\infty.
    \label{eq:price-residual}
\end{equation}
The second measures the largest remaining unilateral profit improvement,
\begin{equation}
    r_\pi^{(n)}
    =
    \max_{s\in A}
    \left[
        \pi_s\!\left(
            BR_s(p_{-s}^{(n)};A),
            p_{-s}^{(n)};
            A
        \right)
        -
        \pi_s(p^{(n)};A)
    \right]_+.
    \label{eq:profit-residual}
\end{equation}
A sweep is accepted when
\begin{equation}
    r_p^{(n)}\leq\varepsilon_p,
    \qquad
    r_\pi^{(n)}\leq\varepsilon_{\pi,p},
    \label{eq:price-convergence}
\end{equation}
with $\varepsilon_p=10^{-5}$ and $\varepsilon_{\pi,p}=10^{-4}$, and all capacity constraints are satisfied. One further sweep then recomputes every best response on a grid of 1,200 prices. The equilibrium is accepted only if this sweep also satisfies \eqref{eq:price-convergence}, so $r_\pi$ is a certificate that no operator can gain more than $\varepsilon_{\pi,p}$ by any unilateral price change, up to the grid resolution and the Newton refinement.

\paragraph{Multiple price equilibria.}
For each reported service set, the price stage is solved again from six further initial price vectors: prices just above cost, half-way between cost and the equilibrium, 1.5 times the equilibrium prices, and three random vectors. Equilibria whose prices differ by more than 0.5 are counted as distinct (\cref{app:diagnostics}).

\subsection{Search over market structures}
\label{sec:search}

The number of possible active service sets grows exponentially with the size of the candidate set. The market-structure problem is therefore treated as a discrete neighbourhood search around repeatedly solved price equilibria.

Starting from an active set $A$, the algorithm first solves $p^*(A)$. It then considers three classes of unilateral structural moves, in this order.

\paragraph{Exit.}
If one or more active services have negative equilibrium profit, the service with the lowest profit is removed and the price equilibrium is solved again.

\paragraph{Route replacement.}
For each incumbent $s\in A$, candidate replacements $r\notin A$ are evaluated by solving the post-replacement price equilibrium for $A^{s\rightarrow r}$. The gain is given by \eqref{eq:switch-gain}. If several profitable replacements exist, a best-improvement rule implements the move with the largest gain.

\paragraph{Entry.}
For each candidate entrant $r\notin A$, the post-entry price equilibrium for $A^{+r}$ is solved and the entrant profit is recorded. If several profitable entrants exist, a best-improvement rule admits the entrant with the largest positive equilibrium profit.

Within each class, candidates are evaluated in decreasing order of an estimate of their gain: the deviator's best response to the current prices of the other services, capped by the bound $B_r$ in \eqref{eq:standalone-bound}. Evaluation stops when the estimate of the next candidate cannot beat the best move already found. The estimate ignores the price response of the incumbents, so the search can overlook a profitable move. The verification step below does not use it.

After any accepted move, the price equilibrium is solved again and the neighbourhood search restarts. This ensures that each discrete decision is evaluated against a strategically consistent market outcome. The search records every active set it visits. A return to a previously visited set is reported as a cycle rather than broken by an arbitrary rule.

\subsection{Candidate set and exact pruning}
\label{sec:pruning}

The candidate set $S$ can be constructed in different ways. In small networks, cyclic rotations up to a chosen number of port calls can be enumerated and combined with vessel classes and speed levels; candidates that cannot carry any OD market within its transit-time limit are dropped. In larger networks, route candidates can be generated selectively using demand-guided search, mixed-integer route synthesis, or column-generation-like pricing problems. In all cases, the equilibrium concept remains conditional on the resulting set $S$. Route generation by mixed-integer synthesis is a heuristic unless it is derived from a reduced-cost pricing problem with an exact master formulation, and the results reported here do not rely on it.

Full equilibrium evaluation of every possible entry or replacement can be expensive, because each candidate requires another price-equilibrium solve. The bound \eqref{eq:standalone-bound} makes most of these solves unnecessary without any loss of rigour. Entry by $r$ cannot be profitable if $B_r\leq\varepsilon_\pi$, and replacement of $s$ by $r$ cannot be profitable if $B_r-\pi_s^*(A)\leq\varepsilon_\pi$, whatever the rest of the market structure. These deviations are skipped; all others are solved. Unlike a screening score, the bound cannot discard a profitable deviation, so the equilibrium certificate does not depend on it.

\subsection{Independent equilibrium verification}
\label{sec:verification}

A local search may terminate because no profitable move is found in the order in which candidates were evaluated, even though a profitable deviation exists elsewhere in the candidate set. The model therefore separates \emph{search} from \emph{verification}.

After the outer search has reached a candidate endpoint, a verification pass tests the conditions \eqref{eq:eq-exit}--\eqref{eq:eq-switch} for every active service, every inactive candidate and every incumbent--candidate pair, skipping only the deviations excluded by the bound in \cref{sec:pruning}. Clone entry is evaluated as well and reported as a diagnostic. A failed price solve in any deviation is counted, and an endpoint with failed solves is not reported as verified. \Cref{app:diagnostics} reports the numbers of deviations solved and skipped.

If verification identifies a profitable deviation, the most profitable one is implemented, the price equilibrium is recomputed, and the market-structure search resumes. The procedure therefore follows the general logic
\begin{equation}
    \text{solve}
    \;\longrightarrow\;
    \text{verify}
    \;\longrightarrow\;
    \text{repair if necessary}
    \;\longrightarrow\;
    \text{solve again}.
\end{equation}
The final competitive outcome is reported as a verified candidate-set equilibrium once the price certificate, the capacity constraints, and the unilateral service-deviation conditions with $\varepsilon_\pi=1$ are all satisfied.
\Cref{alg:competition} in \cref{app:algorithm} summarises the complete procedure.

\subsection{Multiple starting points and equilibrium multiplicity}

Discrete entry and route decisions can produce several locally stable market structures. The solution therefore does not rely on a single initial active set. The complete market-structure problem is solved from 110 initial configurations: the empty market; the five single services with the largest bounds $B_r$; 100 random sets of one to four services drawn from those with $B_r>\varepsilon_\pi$; and the service sets selected by the monopoly, Ramsey and planner regimes.

Two solutions are classified as the same equilibrium if they contain the same active service set and their equilibrium prices differ by less than 0.5 in absolute or $10^{-4}$ in relative terms. All distinct verified equilibria are reported. The tables report the equilibrium reached from the empty market, which corresponds to sequential entry into an initially unserved market; no selection by welfare or profit is made.

Multistart analysis serves both a computational and an economic purpose. Computationally, it reduces the risk that a local-search artefact is mistaken for the market outcome. Economically, it reveals genuine equilibrium multiplicity arising from fixed costs, differentiated services, and discrete market participation.

\subsection{Monopoly, Ramsey, and social-planner service selection}
\label{sec:planner-search}

The alternative regimes share the same discrete service-set problem but replace the competitive price game by a continuous optimisation problem conditional on $A$.

Let
\begin{equation}
    J^R(A)
\end{equation}
denote the optimised objective value for regime $R$ after solving the continuous problem for service set $A$. For monopoly, $J^M(A)$ is optimised total profit. For Ramsey and the social planner, the objective is optimised welfare subject to the corresponding financial constraints. For a fixed set, the problem is a smooth nonlinear programme in the prices. It is solved by sequential quadratic programming with analytic derivatives from several starting points, including the competitive prices where they exist, and the solution is refined by Newton steps on the Karush--Kuhn--Tucker conditions of the active constraints. Capacity shadow prices and Ramsey budget multipliers are recovered by non-negative least squares on the active constraints, so they do not depend on the solver.

The service-set search uses the full neighbourhood
\begin{equation}
    \mathcal N(A)
    =
    \mathcal N^{\mathrm{add}}(A)
    \cup
    \mathcal N^{\mathrm{remove}}(A)
    \cup
    \mathcal N^{\mathrm{swap}}(A),
\end{equation}
in which an add move activates one inactive service, a remove move deletes one active service, and a swap simultaneously removes one active service and adds one inactive service. A best-improvement rule selects
\begin{equation}
    A^{\mathrm{new}}
    =
    \arg\max_{A'\in\mathcal N(A)}
    J^R(A'),
\end{equation}
provided that
\begin{equation}
    J^R(A^{\mathrm{new}})
    >
    J^R(A)+\varepsilon_J,
\end{equation}
with $\varepsilon_J=1$. Neighbours are evaluated in decreasing order of the bounds \eqref{eq:bound-social}--\eqref{eq:bound-monopoly} and the Ramsey bound, and evaluation stops when no remaining bound exceeds the best objective found by more than $\varepsilon_J$. Each step is therefore an exact best-improvement step over the whole neighbourhood, given the fixed-set solutions. The search terminates when no neighbour improves the objective.

The search starts from the empty set, the four single services with the largest bounds, every verified competitive equilibrium, and the solutions of the regimes solved before it (monopoly, then service-level Ramsey, system-wide Ramsey and the planner). Because each regime is started from the solutions of the regimes it must weakly dominate in \eqref{eq:welfare-ordering-property}, a violation of \cref{prop:ordering} in the results would reveal a failure of the fixed-set optimisation rather than an early stop of the search. The model itself defines the monopoly, Ramsey and planner problems globally over $A\subseteq S$; the search is locally optimal with respect to the full add--remove--swap neighbourhood, and \cref{sec:validation} compares it with full enumeration on reduced instances.

\subsection{Numerical validity and robustness}

Numerical validity is assessed separately for competitive and centrally optimised market regimes. A valid competitive equilibrium must satisfy all active-service capacity constraints and meet the price-equilibrium tolerances in \eqref{eq:price-convergence} on the certificate grid. In addition, no active operator may profitably exit, no inactive service may profitably enter, and no incumbent may profitably replace its service. All entry and replacement deviations must either be evaluated explicitly or excluded by the analytical profit bound.

For monopoly, Ramsey pricing and social planning, numerical validity requires feasible prices and capacities, satisfaction of any applicable cost-recovery constraints, and convergence of the continuous optimisation within the stated tolerances. The selected service set must also be locally optimal with respect to the full add--remove--swap neighbourhood.

Robustness is assessed through sensitivity to numerical tolerances, candidate set construction, and initial service configurations. These checks distinguish numerical reliability from the underlying economic assumptions, which is particularly important when combining nonlinear continuous optimisation with discrete market-structure decisions.

\newpage
\section{Computational study}
\label{sec:study}

This section tests the model on a small synthetic network of Baltic and North Sea ports. The aim is not to describe an actual market. The study checks that the solution method works, compares the five market regimes under identical conditions, and examines how the outcomes respond to changes in fuel prices, demand and vessel capacity. All numbers are produced directly by the model code. The main results use multinomial logit demand (\cref{sec:validation,sec:base,sec:scenarios,sec:robustness}), and a second set of runs repeats the scenarios with nested logit demand (\cref{sec:nested}).

\subsection{Sample network}

The sample network is modelled on the Baltic instance of the LINER-LIB benchmark suite \parencite{brouer2014}, a regional feeder network that is widely used to test methods for liner-shipping network design. It keeps the main features of that instance, namely cyclic services, two vessel classes of similar size and directed origin--destination demand, but it is reduced to six ports and ten markets. The ports are Stockholm, Helsinki, Tallinn, Riga, Lübeck and Rotterdam, so the network combines short Baltic routes with longer routes to the German and Dutch gateway ports. Unlike LINER-LIB, the network includes shipper choice. Freight belongs to one of two groups, time-sensitive higher-value goods and less time-sensitive goods, each with its own logit parameters, and every market has an outside option representing road, rail or other alternatives. Demand differs by direction, and the outside option is more attractive in markets where land alternatives are assumed to be stronger. The network is therefore not a solution to the original benchmark but a reduced variant, designed so that every equilibrium can be verified and the regimes can be compared under controlled conditions.

Two vessel classes and three sailing speeds are available. Candidate services are formed by listing every loop that calls at two to four ports, in either direction, which gives 145 loops. Combining each loop with the two vessel classes and three speeds gives 870 services. Of these, \fmBaseCandidates{} can carry at least one origin--destination pair within its maximum transit time, and they form the candidate set $S$. All parameter values, including the utility normalisation, are listed in \cref{app:parameters}.

The parameter values are synthetic but chosen to make the mechanisms of interest visible: freight groups that value price and time differently, fixed service costs, strategic pricing, entry and exit, the trade-off between speed and fuel cost, and limited vessel capacity. Prices, welfare levels and modal shares are therefore not estimates for any real market. What the study can show is whether the model responds sensibly to changes in market conditions and whether the regimes produce outcomes that can be interpreted.

\subsection{Scenarios}

To examine how market regimes respond to changes in economic and operating conditions, four scenarios are considered, as summarised in \cref{tab:scenarios}. The base case provides a reference for comparing market regimes, while the remaining scenarios examine the effects of higher fuel costs, increased freight demand and reduced vessel capacity.

\begin{table}[H]
  \centering
  \caption{Scenario definitions.}
  \label{tab:scenarios}
  \begin{tabularx}{\textwidth}{>{\raggedright\arraybackslash}p{0.19\textwidth} X X}
    \toprule
    Scenario & Change relative to the base case & Purpose \\
    \midrule
    Base & Reference demand, costs, and capacity & Compare market regimes \\
    High fuel & Fuel price increased by 50\% & Test operating-cost response \\
    High demand & All OD demands increased by 25\% & Test entry and network expansion \\
    Low capacity & Vessel capacity reduced by 50\% & Test scarcity and market restructuring \\
    \bottomrule
  \end{tabularx}
\end{table}

The same candidate-generation procedure and behavioural parameters are maintained across scenarios, allowing the effects of fuel costs, demand and capacity to be examined separately. In the low-capacity scenario, vessel capacity is halved while vessel and fuel costs remain unchanged, representing a supply disruption rather than a shift to smaller vessels.

\subsection{Implementation and validation}
\label{sec:validation}

The model is implemented in Python 3.12 using NumPy, SciPy and pandas and requires no commercial solver. The one-dimensional best responses and the evaluation of shares and profits are compiled with Numba, the fixed-set problems of the monopoly, Ramsey and planner regimes are solved with SciPy's sequential quadratic programming routine, and the independent price equilibria needed for verification are solved in parallel.
A complete run covers all four scenarios, the five regimes, 110 competitive
starts per scenario and the validation checks; \cref{app:diagnostics} reports
the run times on a 24-core desktop computer.
Every run records its settings, input data and software versions, and all tables and numbers in this section are generated directly from these runs.

The solver was validated in two ways. Where closed-form prices exist (one freight group and slack capacity), the numerical mark-ups match \eqref{eq:markup-competition}--\eqref{eq:markup-monopoly} and marginal-cost pricing to within \fmValCompDev{} per FFE or better, under both multinomial and nested logit. On two small networks where every set of services can be enumerated, the search reaches the true equilibrium and the global optimum of every regime from every starting point. Both small networks have a single equilibrium, so this confirms that equilibria are found and certified correctly, not that the multistart finds every equilibrium when there are several.

\subsection{Base case}
\label{sec:base}

\begin{figure}[H]
  \centering
  \includegraphics[width=\textwidth]{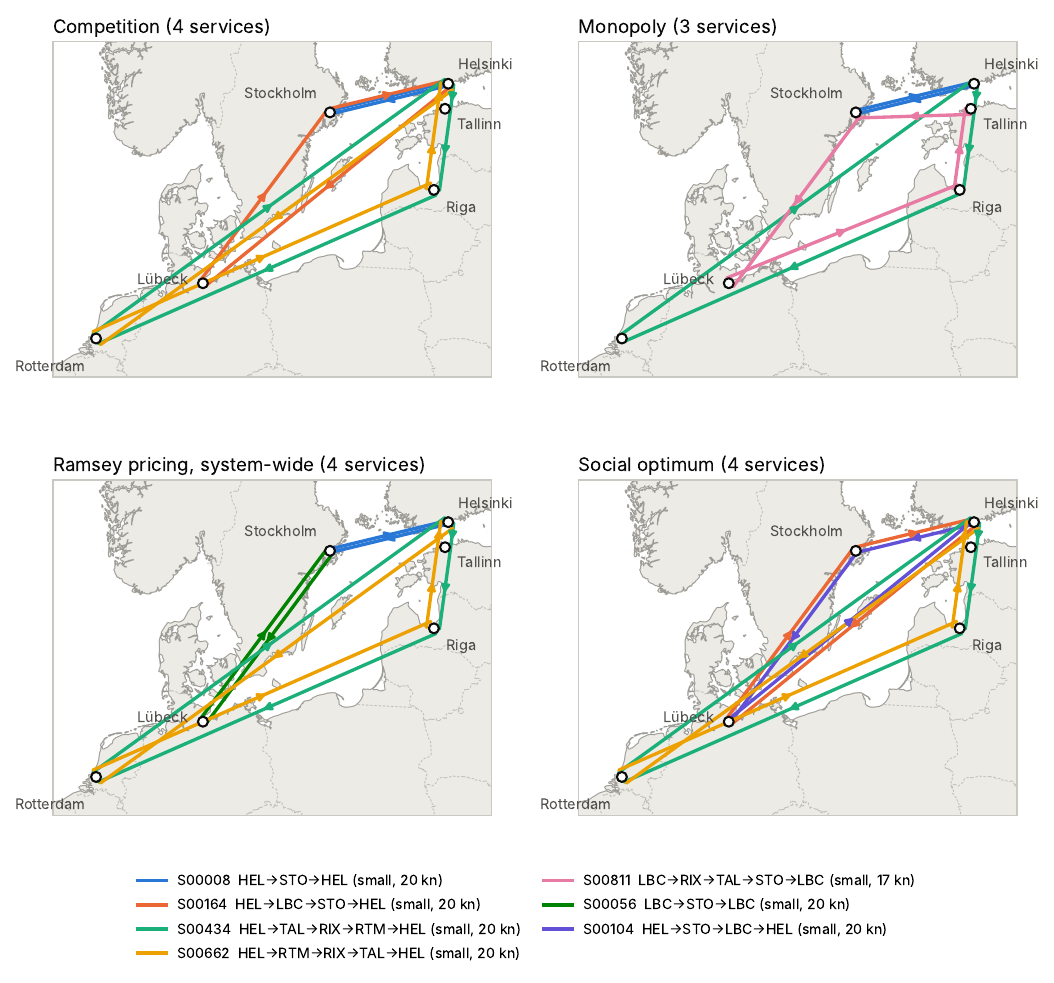}
  \caption{Active liner services in the base case under competition, monopoly, system-wide Ramsey pricing, and the social optimum. Each line is a weekly cyclic service and the arrows give its direction; services that call at the same ports are drawn side by side. A service has the same colour in every panel.}
  \label{fig:active-services}
\end{figure}

\Cref{tab:base-results} reports the base-case outcomes, \cref{tab:base-services} lists the selected services, and \cref{fig:active-services} shows them on a map. Welfare components are reported relative to competition, which makes them invariant to the utility normalisation (\cref{prop:normalisation}). Every selected service in every regime uses the smaller vessel class, whose capacity the demand does not fill.

\begin{table}[H]
  \centering
  \caption{Base-case outcomes by market regime. $\Delta CS$ and $\Delta W$ are differences from competition in model currency units per week; producer surplus is before any subsidy. Average sea price: volume-weighted mean price per FFE over all maritime flows.}
  \label{tab:base-results}
  \small
  \setlength{\tabcolsep}{4pt}
\begin{tabularx}{\textwidth}{>{\raggedright\arraybackslash}p{0.2\textwidth} >{\centering\arraybackslash}p{0.08\textwidth} >{\raggedleft\arraybackslash}X >{\raggedleft\arraybackslash}X >{\raggedleft\arraybackslash}X >{\raggedleft\arraybackslash}X >{\raggedleft\arraybackslash}X }
    \toprule
    Regime & Services & Sea share (\%) & Avg.\ sea price & Producer surplus & $\Delta CS$ & $\Delta W$ \\
    \midrule
    Competition & 4 & 58.7 & 896 & \num{190827} & 0 & 0 \\
    Monopoly & 3 & 45.5 & 1019 & \num{336262} & \num{-353780} & \num{-208345} \\
    Ramsey: service-level & 4 & 70.8 & 635 & $\approx 0$ & \num{401563} & \num{210736} \\
    Ramsey: system-wide & 4 & 73.7 & 594 & $\approx 0$ & \num{435706} & \num{244878} \\
    Social optimum & 4 & 89.9 & 56 & \num{-991605} & \num{1533542} & \num{351111} \\
    \bottomrule
\end{tabularx}
 \end{table}

\begin{table}[H]
  \centering
  \caption{Services selected in the base case (identifier: rotation/vessel class/speed). A rotation returns to its first port.}
  \label{tab:base-services}
  \footnotesize
\begin{tabularx}{\textwidth}{>{\raggedright\arraybackslash}p{0.24\textwidth} X}
    \toprule
    Regime & Services (id: rotation/vessel/speed) \\
    \midrule
    Competition & S00008:HEL-STO/SMALL/20kn;\newline S00164:HEL-LBC-STO/SMALL/20kn;\newline S00434:HEL-TAL-RIX-RTM/SMALL/20kn;\newline S00662:HEL-RTM-RIX-TAL/SMALL/20kn \\
    Monopoly & S00008:HEL-STO/SMALL/20kn;\newline S00434:HEL-TAL-RIX-RTM/SMALL/20kn;\newline S00811:LBC-RIX-TAL-STO/SMALL/17kn \\
    Ramsey: service-level & S00008:HEL-STO/SMALL/20kn;\newline S00164:HEL-LBC-STO/SMALL/20kn;\newline S00434:HEL-TAL-RIX-RTM/SMALL/20kn;\newline S00662:HEL-RTM-RIX-TAL/SMALL/20kn \\
    Ramsey: system-wide & S00008:HEL-STO/SMALL/20kn;\newline S00056:LBC-STO/SMALL/20kn;\newline S00434:HEL-TAL-RIX-RTM/SMALL/20kn;\newline S00662:HEL-RTM-RIX-TAL/SMALL/20kn \\
    Social optimum & S00104:HEL-STO-LBC/SMALL/20kn;\newline S00164:HEL-LBC-STO/SMALL/20kn;\newline S00434:HEL-TAL-RIX-RTM/SMALL/20kn;\newline S00662:HEL-RTM-RIX-TAL/SMALL/20kn \\
    \bottomrule
\end{tabularx}
 \end{table}

The base case satisfies the ordering of \cref{prop:ordering}:
\begin{equation}
  W^{\SO}>W^{R,\mathrm{sys}}>W^{R,\mathrm{serv}}>W^{C}>W^{M}.
  \label{eq:welfare-ordering}
\end{equation}
Only the position of competition relative to monopoly is informative here; the
rest holds by construction and, since each search starts from the solutions it
must dominate (\cref{sec:planner-search}), checks the fixed-set solvers rather
than the service-set searches, which are checked by full enumeration on reduced
instances (\cref{sec:validation}).

Competition sustains \fmBaseCompN{} services and monopoly \fmBaseMonoN. The difference lies in overlap. In the competitive equilibrium, \fmBaseCompContested{} of the 10 OD markets served by sea are served by two or more services. The monopolist, which internalises the business that one of its services takes from another, chooses services that never serve the same market. The volume-weighted average maritime price is \fmBaseCompPrice{} under competition and \fmBaseMonoPrice{} under monopoly, \fmBasePriceGapPct{} per cent lower under competition, and the maritime share is \fmBaseCompShare{} against \fmBaseMonoShare{} per cent. Relative to competition, monopoly changes consumer surplus by \fmBaseMonodCS{} and welfare by \fmBaseMonodW{}; its producer surplus is \fmBaseMonoPS{} against \fmBaseCompPS{} under competition.

Ramsey pricing gains about as much again. Service-level Ramsey pricing selects the competitive service set itself, so its welfare gain of \fmBaseRamServdW{} over competition is a pure pricing effect: each service charges the lowest prices that still cover its own costs. System-wide Ramsey pricing replaces one competitive service by a direct Lübeck--Stockholm service and lets the Helsinki--Stockholm service, which earns a surplus, cross-subsidise the other three, which run deficits at Ramsey prices. It lowers the average price to \fmBaseRamSysPrice, raises the maritime share to \fmBaseRamSysShare{} per cent, and raises welfare by \fmBaseRamSysdW{} over competition. That is \fmBaseRamseyShareOfPlannerGain{} per cent of the planner's gain of \fmBaseSocialdW, achieved without a deficit. The two Ramsey variants differ by \fmBaseRamseyGap, far above the search tolerance $\varepsilon_J=1$. Under service-level recovery the budget multipliers of the four services range from 0.13 to 1.43: the welfare cost of making each service break even differs by a factor of eleven across services, and pooling the constraint equalises it at 0.41. The planner prices at marginal cost (average price \fmBaseSocialPrice), carries \fmBaseSocialShare{} per cent of the freight by sea, and runs a producer surplus of \fmBaseSocialPS{}, the fixed costs that prices do not cover.

\Cref{prop:fuel,prop:demand,prop:normalisation} are confirmed numerically. With the base-case service sets held fixed, the 50 per cent fuel-price increase leaves competitive, monopoly and planner prices unchanged and lowers competitive welfare by exactly the increase in fixed cost, \fmFuelFixedSetWelfareLoss. The 25 per cent demand increase leaves the same prices unchanged, because no leg binds. Adding one to every utility leaves all prices unchanged to within $10^{-8}$ and shifts consumer surplus by $\sum_{k,g}D_{kg}/\beta_g$ (\num{1075000}) in every regime.

\subsection{Scenario results}
\label{sec:scenarios}

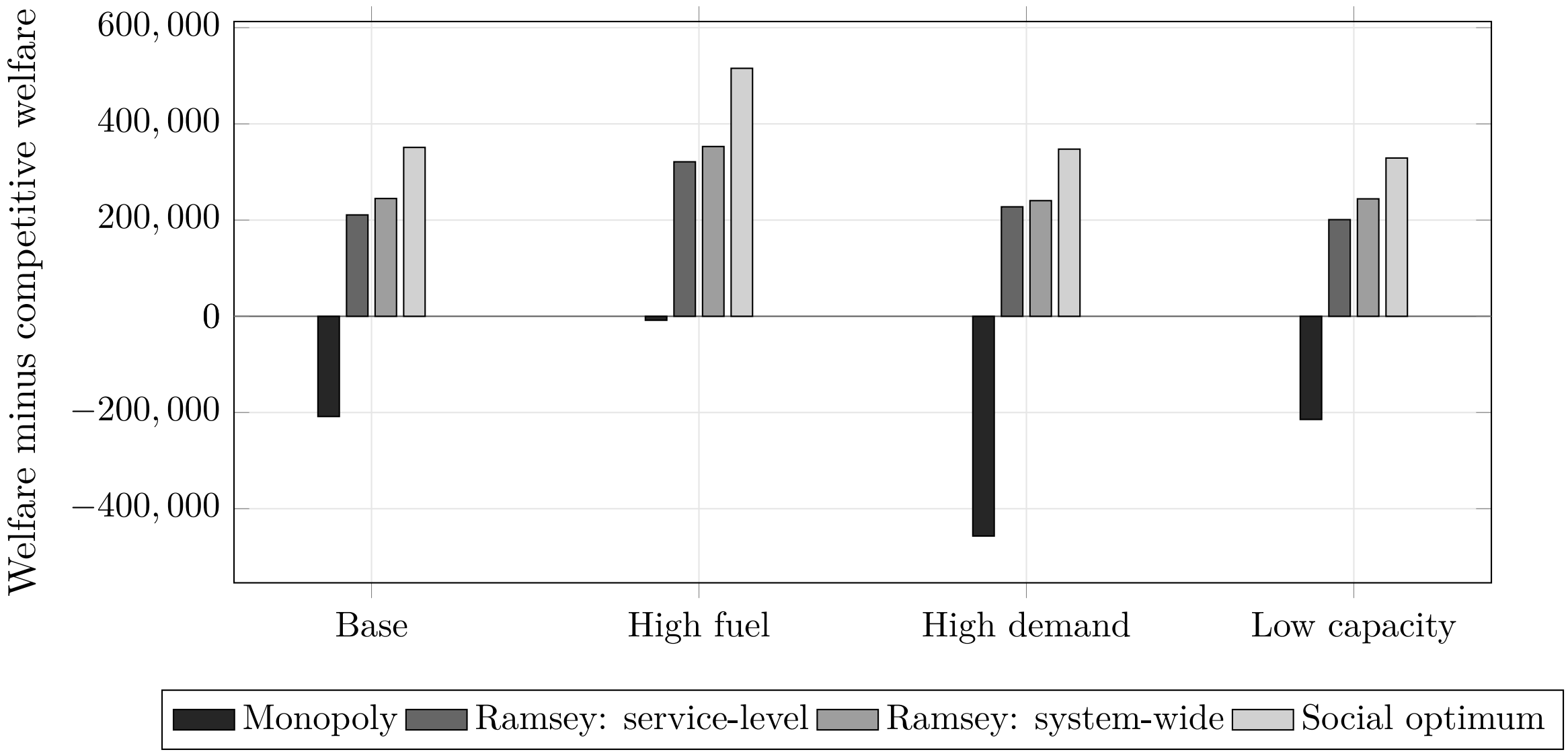
\begin{figure}[H]
  \centering
  \begin{tikzpicture}
  \begin{axis}[
    ybar,
    area legend,
    cycle list name=greys,
    width=0.93\textwidth,
    height=7.1cm,
    bar width=6pt,
    ylabel={Welfare minus competitive welfare},
    symbolic x coords={Base,High fuel,High demand,Low capacity},
    xtick=data,
    x tick label style={font=\small},
    y tick label style={/pgf/number format/fixed,/pgf/number format/1000 sep={{,}},font=\small},
    scaled y ticks=false,
    legend style={at={(0.5,-0.19)},anchor=north,legend columns=4,font=\small},
    enlarge x limits=0.14,
    grid=major,
    grid style={draw=black!10},
    extra y ticks={0},
    extra y tick style={grid style={draw=black!60}},
  ]
\addplot coordinates {(Base,-208345) (High fuel,-8265) (High demand,-456840) (Low capacity,-214477)}; 
\addplot coordinates {(Base,210736) (High fuel,321155) (High demand,227566) (Low capacity,200934)}; 
\addplot coordinates {(Base,244878) (High fuel,352939) (High demand,240419) (Low capacity,244169)}; 
\addplot coordinates {(Base,351111) (High fuel,515532) (High demand,347525) (Low capacity,329036)}; 
     \legend{Monopoly,Ramsey: service-level,Ramsey: system-wide,Social optimum}
  \end{axis}
  \end{tikzpicture}
  \caption{Welfare of each regime minus competitive welfare in the same scenario (model currency units per week). Competition is the equilibrium reached from the empty market.}
  \label{fig:welfare-gap}
\end{figure}

\begin{figure}[H]
  \centering
  \begin{tikzpicture}
  \begin{axis}[
    ybar,
    area legend,
    cycle list name=greys,
    width=0.93\textwidth,
    height=7.1cm,
    bar width=5pt,
    ylabel={Maritime market share (\%)},
    ymin=0,ymax=100,
    symbolic x coords={Base,High fuel,High demand,Low capacity},
    xtick=data,
    x tick label style={font=\small},
    y tick label style={font=\small},
    legend style={at={(0.5,-0.19)},anchor=north,legend columns=3,font=\small},
    enlarge x limits=0.14,
    grid=major,
    grid style={draw=black!10},
  ]
\addplot coordinates {(Base,58.7) (High fuel,45.4) (High demand,66.0) (Low capacity,58.7)}; 
\addplot coordinates {(Base,45.5) (High fuel,45.4) (High demand,45.8) (Low capacity,45.1)}; 
\addplot coordinates {(Base,70.8) (High fuel,66.8) (High demand,77.5) (Low capacity,70.1)}; 
\addplot coordinates {(Base,73.7) (High fuel,68.0) (High demand,78.7) (Low capacity,73.6)}; 
\addplot coordinates {(Base,89.9) (High fuel,89.1) (High demand,91.1) (Low capacity,90.4)}; 
     \legend{Competition,Monopoly,Ramsey: service-level,Ramsey: system-wide,Social optimum}
  \end{axis}
  \end{tikzpicture}
  \caption{Maritime market share across market regimes and scenarios.}
  \label{fig:sea-share-scenarios}
\end{figure}
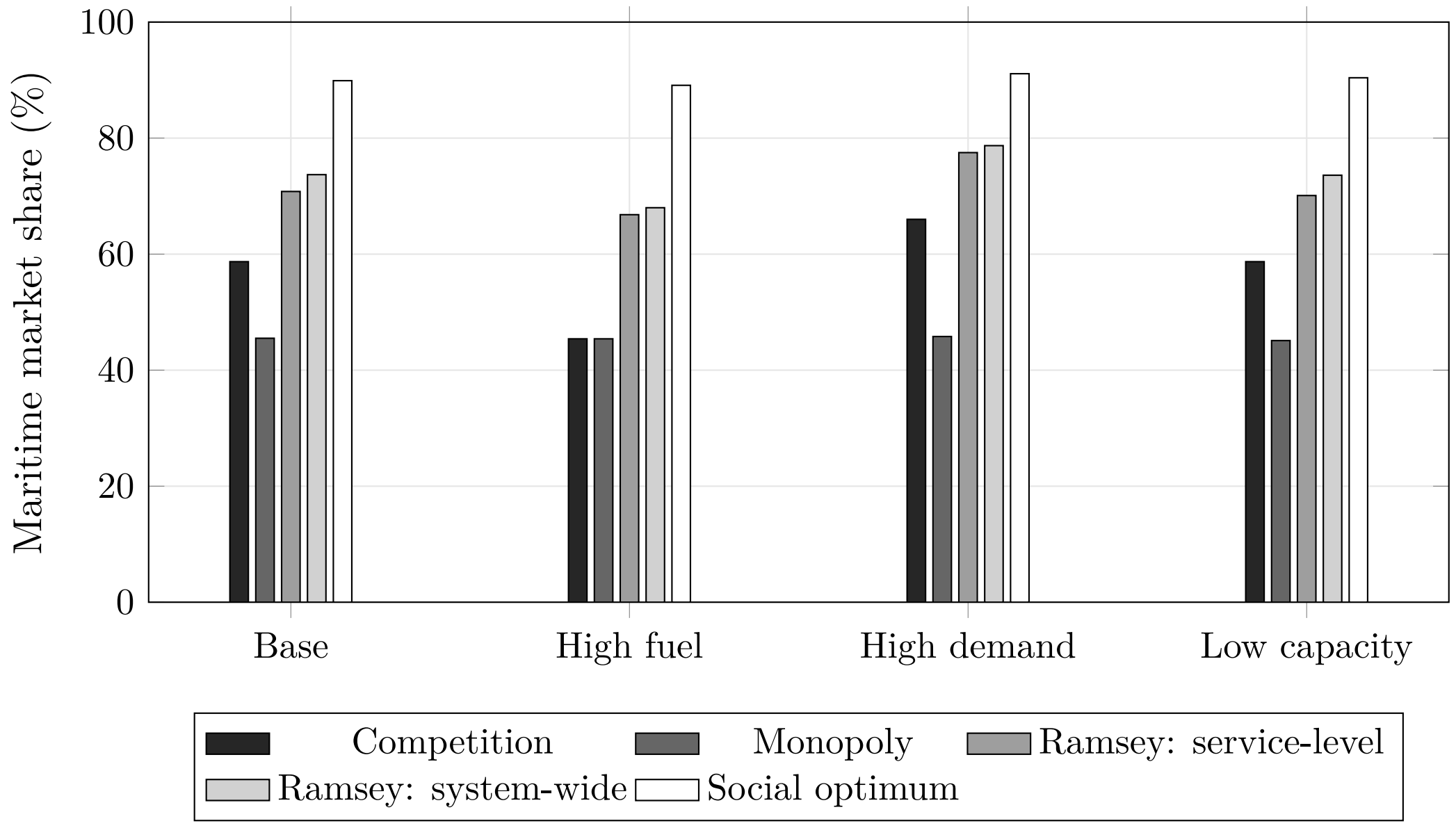

\Cref{fig:welfare-gap} shows, for each scenario, the welfare of each regime relative to competition in the same scenario; these differences hold demand fixed and are therefore invariant to the utility normalisation. \Cref{fig:sea-share-scenarios} reports maritime shares, \cref{tab:welfare-changes} reports the change in welfare relative to the base case, and \cref{tab:competition-structure} summarises the competitive outcomes.

\begin{table}[H]
  \centering
  \caption{Change in welfare relative to the base case (model currency units per week). The fuel and capacity scenarios hold demand fixed, so their rows are invariant to the utility normalisation; the high-demand row is not and includes the mechanical scale effect of \cref{prop:demand}.}
  \label{tab:welfare-changes}
  \small
\begin{tabularx}{\textwidth}{>{\raggedright\arraybackslash}p{0.16\textwidth} >{\raggedleft\arraybackslash}X >{\raggedleft\arraybackslash}X >{\raggedleft\arraybackslash}X >{\raggedleft\arraybackslash}X >{\raggedleft\arraybackslash}X }
    \toprule
    Scenario & Competition & Monopoly & Ramsey: service-level & Ramsey: system-wide & Social optimum \\
    \midrule
    High fuel & \num{-274388} & \num{-74308} & \num{-163968} & \num{-166327} & \num{-109966} \\
    High demand & \num{1110290} & \num{861795} & \num{1127120} & \num{1105831} & \num{1106704} \\
    Low capacity & 0 & \num{-6132} & \num{-9802} & \num{-709} & \num{-22074} \\
    \bottomrule
\end{tabularx}
 \end{table}

\begin{table}[H]
  \centering
  \caption{Competitive outcomes by scenario. The first three columns describe the equilibrium reached from the empty market: number of services, OD markets served by two or more services (of the markets served by sea), and its welfare advantage over monopoly. The last three columns cover all distinct verified equilibria found from the 110 starts.}
  \label{tab:competition-structure}
  \footnotesize
\begin{tabularx}{\textwidth}{>{\raggedright\arraybackslash}p{0.14\textwidth} >{\raggedleft\arraybackslash}X >{\raggedleft\arraybackslash}X >{\raggedleft\arraybackslash}X >{\raggedleft\arraybackslash}X >{\raggedleft\arraybackslash}X >{\raggedleft\arraybackslash}p{0.22\textwidth}}
    \toprule
    Scenario & Services & Contested markets & $W^C-W^M$ & Equilibria (starts) & Equilibria with a contested market & Welfare range of equilibria \\
    \midrule
    Base & 4 & 5/10 & \num{208345} & 1 (110) & 1 & -- \\
    High fuel & 3 & 1/9 & \num{8265} & 2 (110) & 2 & \num{2659387}--\num{2722596} \\
    High demand & 5 & 8/10 & \num{456840} & 1 (110) & 1 & -- \\
    Low capacity & 4 & 5/10 & \num{214477} & 1 (110) & 1 & -- \\
    \bottomrule
\end{tabularx}
 \end{table}

\subsubsection{Fuel-cost shock}

A 50 per cent fuel-price increase lowers welfare in every regime. By \cref{prop:fuel}, competitive and monopoly prices for a given set of services do not depend on the fuel price, so every price change in these regimes is a service-set effect, and the competitive set does change. At the higher fuel price the base-case set is no longer an equilibrium, because one of its four services makes a loss. In the equilibrium reached from the empty market, that service has left, and the Helsinki--Stockholm and Helsinki--Lübeck--Stockholm services sail at 17 instead of 20 knots. Competition now runs \fmHighFuelCompN{} services, the maritime share falls from \fmBaseCompShare{} to \fmHighFuelCompShare{} per cent, and the average price rises from \fmBaseCompPrice{} to \fmHighFuelCompPrice.

The structural response amplifies the welfare loss rather than cushioning it. With the base-case services held fixed, welfare would fall by the increase in fixed cost, \fmFuelFixedSetWelfareLoss. With the services adjusted, competitive welfare falls by \fmHighFuelCompdWbase{} (\cref{tab:welfare-changes}). Exit and slow steaming save the operators fuel, but shippers lose a service and face longer transit times, and these losses are larger than the fuel saved. Monopoly makes only one change, slowing its Helsinki--Stockholm service, and loses \fmHighFuelMonodWbase. The welfare difference between monopoly and competition, $W^M-W^C$, therefore narrows from \fmBaseMonodW{} to \fmHighFuelMonodW.

The fuel scenario also has a second competitive equilibrium, with four services, slower vessels and higher welfare (\cref{tab:competition-structure}): competitive welfare is \fmHighFuelCompWmin{} in the equilibrium reached from the empty market and \fmHighFuelCompWmax{} in the other. The size of the welfare loss from the fuel shock therefore depends on which equilibrium the market coordinates on. Ramsey prices respond directly, because fuel enters the cost-recovery constraint. The planner keeps prices at marginal cost, changes two of its four services, and its maritime share barely moves (\fmBaseSocialShare{} to \fmHighFuelSocialShare{} per cent); the higher fuel bill is financed through transfers.

\subsubsection{Demand increase}

A 25 per cent increase in all OD demands produces the clearest entry effect. Competition grows from \fmBaseCompN{} to \fmHighDemandCompN{} services: a Helsinki--Stockholm--Lübeck service at 17 knots enters, running the reverse of an incumbent rotation, and eight of the ten markets served are now contested. The maritime share rises from \fmBaseCompShare{} to \fmHighDemandCompShare{} per cent and the average maritime price falls from \fmBaseCompPrice{} to \fmHighDemandCompPrice. By \cref{prop:demand}, uniform demand growth does not change prices for a given set of services when capacity is slack, and the model confirms that the base-case prices are unchanged at the higher demand. The price fall is therefore an entry effect: the larger market covers the fixed cost of a fifth service, which then competes with the incumbents.

Monopoly keeps three services but replaces its Lübeck--Riga--Tallinn--Stockholm service at 17 knots by the same rotation at 20 knots. Its average price rises slightly, to \fmHighDemandMonoPrice, which again is a service-set effect, and its maritime share barely changes. The welfare difference $W^M-W^C$ widens from \fmBaseMonodW{} to \fmHighDemandMonodW. Ramsey pricing and the planner also expand to five services. Welfare rises in every regime, but these increases include the mechanical scale effect and depend on the normalisation. Within the scenario, the advantage of system-wide Ramsey pricing over competition narrows slightly, from \fmBaseRamSysdW{} to \fmHighDemandRamSysdW, because competition itself has become more intense.

\subsubsection{Capacity reduction}

Halving vessel capacity leaves the competitive equilibrium unchanged: the same \fmLowCapacityCompN{} services, the same prices and the same welfare. No leg of the base-case competitive network was more than half full, so the base-case prices remain feasible; the largest leg load factor doubles from \fmBaseCompLF{} to \fmLowCapacityCompLF, and no leg binds. Verification confirms that no entry or replacement becomes profitable, and no other equilibrium is found from the 110 starts.

The other regimes do hit the constraint. The monopolist concentrates traffic on three services, fills one leg, and raises its average price to \fmLowCapacityMonoPrice; it loses \fmLowCapacityMonodWbase. Service-level Ramsey pricing loses \fmLowCapacityRamServdWbase{} and system-wide Ramsey pricing \fmLowCapacityRamSysdWbase, so the value of pooling the budget constraint rises when capacity is scarce. The planner, whose low prices generate the most traffic, has two binding legs and responds by rearranging its network into five services; it loses \fmLowCapacitySocialdWbase. Capacity scarcity therefore matters for the regimes that load their services most heavily: the monopolist, with fewer services, and the planner, with the lowest prices.
In this instance, competitive loads stay below capacity under multinomial logit.

\subsection{Equilibrium robustness}
\label{sec:robustness}

The competitive solution was computed from 110 initial market structures in each scenario. Every start ended at a verified equilibrium, and no search cycled. In the base, high-demand and low-capacity scenarios all starts ended at the same equilibrium; in the high-fuel scenario two distinct equilibria were found (\cref{tab:competition-structure}). In every reported equilibrium the most profitable unilateral deviation loses money: its gain lies between \fmHighFuelMaxGain{} and \fmBaseMaxGain{} per week, against a tolerance of $\varepsilon_\pi=1$; clone entry would also make a loss; and the price stage has only one equilibrium among the seven initial price vectors tried. \Cref{app:diagnostics} reports the candidate-set size, the deviations solved and excluded by the bound, residuals, and run times.
Verification is conditional on the candidate set: a larger or differently generated pool can contain profitable routes that the current pool does not.

\subsection{Nested-logit demand}
\label{sec:nested}

\begin{figure}[H]
  \centering
  \begin{tikzpicture}
  \begin{axis}[
    width=0.93\textwidth,
    height=7cm,
    xlabel={Nest parameter $\rho$},
    ylabel={Welfare minus competitive welfare},
    x dir=reverse,
    xtick={1,0.9,0.8,0.7,0.6,0.5,0.4},
    xticklabels={1.0,0.9,0.8,0.7,0.6,0.5,0.4},
    y tick label style={/pgf/number format/fixed,/pgf/number format/1000 sep={\mathord{,}},font=\small},
    x tick label style={font=\small},
    scaled y ticks=false,
    grid=major,
    grid style={draw=black!10},
    extra y ticks={0},
    extra y tick style={grid style={draw=black!60}},
    legend style={at={(0.5,-0.22)},anchor=north,legend columns=2,font=\small,draw=none},
    legend cell align=left,
    cycle list={
      {black,mark=o,thick},
      {black,mark=square,dashed,thick},
      {black,mark=triangle,densely dotted,thick},
      {black,mark=*,thick},
    },
  ]
\addplot coordinates {(1.00,-208345) (0.90,-188411) (0.80,-5009) (0.70,0) (0.60,0) (0.50,0) (0.40,0)}; 
\addplot coordinates {(1.00,210736) (0.90,186472) (0.80,364984) (0.70,369993) (0.60,369993) (0.50,369993) (0.40,369993)}; 
\addplot coordinates {(1.00,244878) (0.90,236494) (0.80,385349) (0.70,390358) (0.60,390358) (0.50,390358) (0.40,390358)}; 
\addplot coordinates {(1.00,351111) (0.90,345967) (0.80,505464) (0.70,487524) (0.60,451105) (0.50,451105) (0.40,451105)}; 
     \legend{Monopoly,Ramsey: service-level,Ramsey: system-wide,Social optimum}
  \end{axis}
  \end{tikzpicture}
  \caption{Welfare of each regime minus competitive welfare in the base case, for nest parameters from 1 (multinomial logit) to 0.4. Competition is the equilibrium reached from the empty market.}
  \label{fig:nest-welfare}
\end{figure}
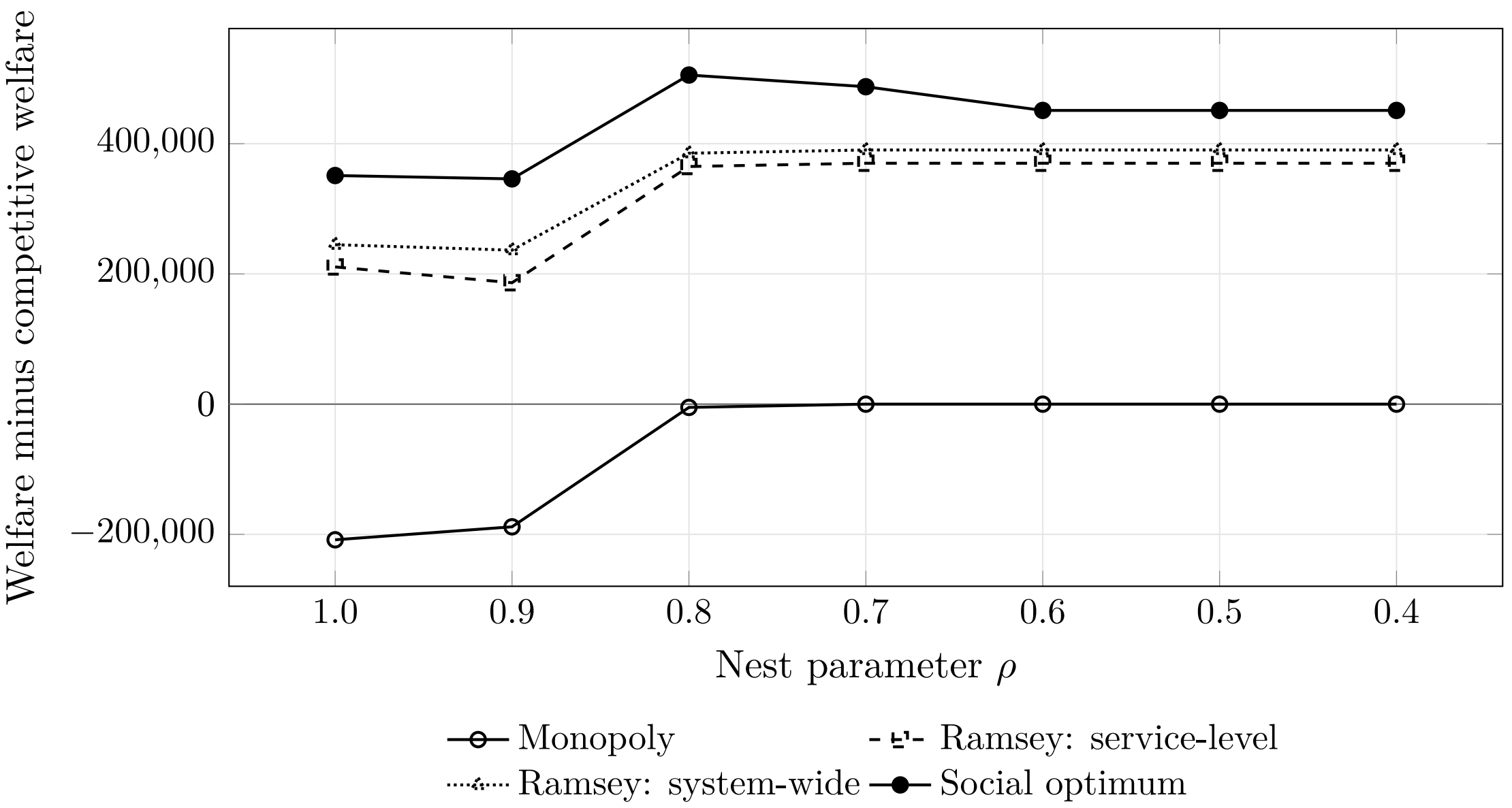

Under multinomial logit, every new service expands the maritime market in proportion to its share, which favours entry into markets that are already served. This section replaces it with the nested logit of \cref{sec:demand}, in which all maritime services belong to one nest. The model is first solved with $\rho=0.6$ for both segments in all four scenarios (\cref{tab:nest-structure}), and then for a range of nest parameters in the base case (\cref{fig:nest-welfare} and \cref{tab:nest-sensitivity}). The value $\rho=0.6$ is illustrative and lies well inside the unit interval, and the range of values shows where the outcome changes. Everything else is unchanged, including the candidate set, the 110 starting configurations and the tolerances.

\begin{table}[H]
  \centering
  \caption{Competitive outcomes under nested-logit demand with $\rho=0.6$. Columns as in \cref{tab:competition-structure}.}
  \label{tab:nest-structure}
  \footnotesize
\begin{tabularx}{\textwidth}{>{\raggedright\arraybackslash}p{0.14\textwidth} >{\raggedleft\arraybackslash}X >{\raggedleft\arraybackslash}X >{\raggedleft\arraybackslash}X >{\raggedleft\arraybackslash}X >{\raggedleft\arraybackslash}X >{\raggedleft\arraybackslash}p{0.22\textwidth}}
    \toprule
    Scenario & Services & Contested markets & $W^C-W^M$ & Equilibria (starts) & Equilibria with a contested market & Welfare range of equilibria \\
    \midrule
    Base & 3 & 0/10 & 0 & 4 (110) & 0 & \num{2698877}--\num{2725430} \\
    High fuel & 3 & 0/10 & 0 & 4 (110) & 0 & \num{2374204}--\num{2651122} \\
    High demand & 3 & 0/10 & 0 & 8 (110) & 1 & \num{3472270}--\num{3805196} \\
    Low capacity & 3 & 0/10 & 0 & 4 (110) & 0 & \num{2652063}--\num{2719298} \\
    \bottomrule
\end{tabularx}
 \end{table}

With $\rho=0.6$, the equilibrium reached from the empty market has \fmNestBaseCompN{} services and no contested market in every scenario, and it coincides with the monopoly solution: the same services, the same prices and the same welfare, so $W^C-W^M=\fmNestBaseCompMinusMonoW$. In the base case it is the monopoly service set of the multinomial-logit model, and by \cref{prop:uncontested} its welfare is the same as there (\fmNestBaseCompW{} against \fmBaseMonoW). Competition has become a set of local monopolies.

\begin{table}[H]
  \centering
  \caption{Base case for a grid of nest parameters. The first five columns describe the competitive equilibrium reached from the empty market; welfare differences are relative to it; the last column counts the distinct verified competitive equilibria found from 110 starts. $\rho=1$ is multinomial logit.}
  \label{tab:nest-sensitivity}
  \footnotesize
  \setlength{\tabcolsep}{3pt}
\begin{tabularx}{\textwidth}{c >{\raggedleft\arraybackslash}X >{\raggedleft\arraybackslash}X >{\raggedleft\arraybackslash}X >{\raggedleft\arraybackslash}X >{\raggedleft\arraybackslash}X >{\raggedleft\arraybackslash}X >{\raggedleft\arraybackslash}X >{\raggedleft\arraybackslash}X }
    \toprule
    $\rho$ & Services & Contested markets & Sea share (\%) & Avg.\ sea price & $W^M-W^C$ & $W^{R,\mathrm{sys}}-W^C$ & $W^{\SO}-W^C$ & Equilibria \\
    \midrule
    1.0 & 4 & 5 & 58.7 & 896 & \num{-208345} & \num{244878} & \num{351111} & 1 \\
    0.9 & 4 & 5 & 59.1 & 865 & \num{-188411} & \num{236494} & \num{345967} & 4 \\
    0.8 & 3 & 1 & 46.1 & 899 & \num{-5009} & \num{385349} & \num{505464} & 4 \\
    0.7 & 3 & 0 & 45.5 & 1019 & 0 & \num{390358} & \num{487524} & 4 \\
    0.6 & 3 & 0 & 45.5 & 1019 & 0 & \num{390358} & \num{451105} & 4 \\
    0.5 & 3 & 0 & 45.5 & 1019 & 0 & \num{390358} & \num{451105} & 5 \\
    0.4 & 3 & 0 & 45.5 & 1019 & 0 & \num{390358} & \num{451105} & 5 \\
    \bottomrule
\end{tabularx}
 \end{table}

The mechanism is the substitution pattern. With $\rho<1$ a service that enters a market already served takes more of its traffic from the incumbent and less from the outside option, and the price competition that follows is fiercer, because demand within the nest is more price-elastic \eqref{eq:markup-competition}. Entry into a served market therefore pays less than entry into an unserved one, and entrants spread out. Once no market is contested, \cref{prop:uncontested} implies that competitive and monopoly prices coincide for that set. The monopolist, which avoids overlap for its own reasons, chooses the same set. The grid in \cref{tab:nest-sensitivity} shows how abruptly this happens. At $\rho=0.9$ the multinomial-logit equilibrium survives with \fmSensNineCompN{} services and lower prices; at $\rho=0.8$ competition has \fmSensEightCompN{} services and \fmSensEightCompContested{} contested market, and its welfare advantage over monopoly has almost vanished; from $\rho=0.7$ down, the equilibrium reached from the empty market reproduces
the monopoly outcome.

Nesting also multiplies the equilibria. With $\rho=0.6$ the multistart finds \fmNestBaseDistinctEq{}, \fmNestHighFuelDistinctEq{}, \fmNestHighDemandDistinctEq{} and \fmNestLowCapacityDistinctEq{} distinct verified equilibria in the four scenarios, against at most two under multinomial logit. All but one consist of services in separate markets. The exception, in the high-demand scenario, is the base-case competitive set of the multinomial-logit model, with five contested markets and the highest welfare of all the equilibria (\fmNestHighDemandCompWmax). It is reached from only a few starts, so the market can sustain overlapping competition but does not arrive there by sequential entry from an empty market.

The planner and Ramsey regimes also select fewer services under nesting: three in the base case instead of four, because a service that largely duplicates another adds less consumer surplus when the services are close substitutes. Because competition is weaker, the gain from system-wide Ramsey pricing over competition is larger than under multinomial logit, \fmNestBaseRamSysdW{} against \fmBaseRamSysdW, and it recovers \fmNestBaseRamseyShareOfPlannerGain{} per cent of the planner's gain. The welfare advantage of competition over monopoly in the main results thus rests on overlapping services, and whether operators choose to overlap depends on how closely shippers regard the services as substitutes. The nest parameter is therefore one of the first quantities an empirical application would need to estimate.

\newpage
\section{Discussion, limitations, and future research}
\label{sec:discussion}

This section discusses the results of \cref{sec:study}, relates them to the literature and to national freight models, and sets out the main limitations of the model and the next steps in its development.

\subsection{Market structure matters beyond operating cost}

The numerical experiments show why the distinction between transport cost and market price matters. With the same demand, costs and candidate services, competition and monopoly select different networks: competitors enter each other's markets, while the monopolist keeps its services apart to avoid cannibalising its own traffic. Overlap is what lowers competitive prices, and it is also what the nested-logit analysis removes. Ramsey pricing generates a much larger maritime share than either commercial regime while covering costs, and the planner expands maritime transport further only by accepting a large deficit.

These differences would not be visible in a model that assigns freight by minimum logistics cost. Strategic pricing creates feedback between service attractiveness, market share, capacity use, and profitability, and endogenous entry and route choice mean that a cost or demand shock can change the set of services offered, not only the flows over a fixed network. \Cref{sec:properties} sharpens this point: for a given set of services, the fuel price and uniform demand growth leave commercial prices unchanged, so their price effects in the experiments run entirely through the service set. The fuel scenario shows that this channel can dominate: the service-set response more than doubles the welfare loss that the fuel bill alone would cause.

\subsection{Entry, fixed costs, and network coverage}

The scenarios show several sides of the same fixed-cost mechanism. Higher demand covers the fixed cost of an additional service and increases competitive coverage and overlap. Higher fuel costs push operators towards slower vessels and make the marginal service unprofitable. Both are consistent with the free-entry literature: the privately sustainable number and configuration of services need not coincide with the welfare-maximising network \parencite{mankiw1986}. In the fuel scenario the market can also settle in either of two equilibria, so the outcome of a shock can depend on how operators coordinate, not only on fundamentals.

Capacity scarcity operates in two ways. A binding leg raises the shadow value of capacity and enters prices through \eqref{eq:economic-mc}. Lower capacity can also change which routes cover their fixed costs, and so change the market structure before the remaining routes are full.
Under multinomial logit, neither happens under competition in the test instance because competitive loads stay below half of the original capacity. Under nested logit with $\rho = 0.6$, competition reproduces the monopoly outcome, and with halved capacity, it fills one leg.

\subsection{Substitution between services}

The nested-logit results qualify the main comparison. Under multinomial logit, every additional service expands the maritime market, so entry into served markets pays and competition lowers prices. When shippers regard maritime services as close substitutes for each other, entrants avoid each other's markets, and from a nest parameter of 0.7 down the competitive outcome reached from the
empty market is identical to monopoly in the test instance.
The same substitution pattern leads the planner and Ramsey regimes to choose fewer, less overlapping services. The welfare gain from competition is therefore not a robust feature of the model structure; it depends on a demand parameter that would have to be estimated, and on the equilibrium that the market reaches, since nesting also multiplies the equilibria.

\subsection{Interpretation of the welfare benchmarks}

The social optimum is an upper benchmark within the model, not an implementable policy. It treats subsidies as costless transfers, excludes the marginal cost of public funds, and represents the outside option through utility rather than an explicit resource-cost model. Its large maritime share should therefore not be read as a desirable modal shift. The Ramsey solutions are more institutionally plausible because they impose cost recovery, but they still presume coordinated welfare-maximising pricing and service selection.

The service-level Ramsey solution uses the competitive service set, so in the base case its advantage over competition is purely a matter of pricing. The system-wide variant gains more by changing the network and letting one profitable service cross-subsidise three others. The widely different budget multipliers under service-level recovery show where the break-even requirement is most costly, and the gap between the two variants grows when capacity is scarce.

\subsection{Relationship to national freight models}

The model is best viewed as a strategic supply-side complement to stochastic freight-demand modules. A Samgods-style logistics model can generate transport-chain attributes and market demand, while the present framework represents how operators convert costs and quality into prices and service provision. A full integration would require consistent zones, commodity classifications, outside alternatives, capacity feedback, and empirical calibration. The current results do not show that the framework is ready for national forecasting. They show that the economic mechanisms can be represented in a computational model and that the resulting equilibria can be subjected to explicit verification.

\subsection{Limitations and extensions}

The main limitations, and the extensions they suggest, are as follows.

\begin{enumerate}[leftmargin=*,label=\arabic*.]
  \item \textbf{Synthetic calibration.} Demand, costs, and behavioural parameters are illustrative. Demand parameters, including the nest parameter, can be estimated from observed freight rates, service attributes, shipment choices, and modal shares, and service costs and capacities from vessel, port, and schedule data.
  \item \textbf{One service per competitive operator.} A multiproduct firm $i$ owning the services $\mathcal S_i$ would maximise $\Pi_i=\sum_{s\in\mathcal S_i}\pi_s$ over the joint price vector of its services. An explicit ownership matrix would represent internal cannibalisation, network complementarities, mergers, and alliances.
  \item \textbf{Substitution patterns.} The nested logit has a single maritime nest. Nests by rotation, by speed or by operator, mixed logit, or latent classes would allow richer substitution, and reliability, frequency, and schedule convenience could enter utility directly.
  \item \textbf{No transshipment.} Shippers select direct services only. Multi-service chains require transfer times, connection feasibility, and rules for allocating revenue and cost across operators, and would create complementarities between services.
  \item \textbf{No global fleet constraint.} Candidate services include the vessels needed for weekly operation, but operators do not compete for a fixed fleet. With $N_f$ vessels of class $f$, a natural constraint is $\sum_{s\in A:f_s=f}m_s\leq N_f$ for all $f\in F$, after which service frequency and vessel deployment could be made endogenous.
  \item \textbf{Static equilibrium and equilibrium selection.} The model finds a weekly equilibrium without sunk costs or investment, and the tables report the equilibrium reached by sequential entry from an empty market. A dynamic extension would separate operating decisions from investments in vessels and terminals and give a principled selection among equilibria \parencite{ericson1995,aguirregabiria2012}.
  \item \textbf{Exogenous total demand and outside option.} Total OD demand is fixed, and the outside transport system does not respond. A broader welfare analysis would model road, rail and other maritime services explicitly, with their costs, emissions and capacity, and would include the marginal cost of public funds when subsidies are used.
  \item \textbf{Numerical scope.} Monopoly, Ramsey, and planner service sets are obtained from multistart searches, and verification is conditional on the candidate set.
  With two segments and a binding leg, the computed best response may not be globally optimal; among the reported equilibria, this occurs only in the nested-logit low-capacity case, where it is not certified (Section~\ref{sec:price-stage}), and it can also occur in deviations evaluated during verification.
\end{enumerate}

\subsection{Future research}

The model can be extended in several directions. The first is to apply it to larger networks, where the candidate services can no longer be listed in full. This calls for more efficient solution methods, such as heuristic and metaheuristic searches over market structures, faster equilibrium algorithms and commercial solvers such as Gurobi. Candidate services could be generated as needed, for example, through mixed-integer or prize-collecting cycle formulations combined with decomposition and column generation. The bound \eqref{eq:standalone-bound} reduces the search space by excluding every route that cannot earn a positive profit even when operated alone. Explicit verification of the equilibrium would remain necessary, so that no profitable deviation is overlooked.

The second direction is a more realistic representation of transport supply, with fleet constraints, multiproduct operators and transhipment between services. Calibration to observed freight rates, services and freight flows would also be needed to assess the model under realistic market conditions.

The third direction concerns the existence and selection of equilibria. Pure-strategy equilibria are not guaranteed, and it remains open under which conditions on demand, costs and capacities a stable service set exists. Small constructed examples could show when cycles of entry and replacement arise, and mixed-strategy or dynamic formulations could describe the market when no pure equilibrium exists. Where several equilibria coexist, as in some of the scenarios here, a dynamic model of entry and investment would give a better basis for choosing among them than sequential entry from an empty market.

In the longer term, the framework could be integrated with national freight transport models. The model has been developed within the HERMES project, which also studies other extensions of national freight modelling, including loop-based representations of maritime transport and network assignment. Combining these components would allow transport prices, service networks and market structure to respond endogenously to changes in freight demand, infrastructure and transport policy, thereby supporting more comprehensive analyses of freight markets under alternative market and regulatory regimes.

\section{Conclusions}
\label{sec:conclusions}

This paper presents a freight-market model that connects shipper choice, strategic pricing, endogenous entry, and liner-service choice. Shippers choose among differentiated maritime services and an outside option; operators set OD prices and can enter, exit, or replace their service; and competition is compared with monopoly, Ramsey pricing, and a planner on a common demand and cost system.

Four findings stand out in the synthetic experiments. First, under multinomial logit, competition sustains more services than monopoly, with overlapping routes, prices \fmBasePriceGapPct{} per cent lower and higher welfare, while cost-recovering Ramsey prices raise the maritime share from \fmBaseCompShare{} to \fmBaseRamSysShare{} per cent at zero deficit. Second, cost and demand shocks change competitive prices only through the service set, as the model structure implies, and this channel can be large: a fuel-price increase leads to exit and slow steaming that more than double its welfare cost, and a demand increase attracts an additional competitor.
Third, under multinomial logit, halving capacity does not change the competitive
outcome, because competitive loads are low, but capacity then binds under
monopoly, both Ramsey regimes and the planner. Fourth, when maritime services are
close substitutes for each other, entrants avoid each other's markets and the
equilibrium reached by sequential entry reproduces the monopoly outcome, so the
gains from competition rest on a substitution parameter that must be estimated.

The methodological lesson is that convergence of a heuristic search is not evidence of equilibrium. Best responses must be global, and every unilateral entry, exit and route replacement must be tested against the post-deviation price equilibrium, or excluded by a valid bound. The solve--verify--repair procedure makes this operational, and multistart shows that several equilibria can coexist. The paper is a framework and computational study rather than a calibrated policy model. Targeted route generation, fleet constraints, multiproduct operators, transshipment, and empirical calibration are the steps towards applied use.

\section*{Acknowledgements}

This work was carried out within the HERMES project (High Resolution Freight Network Assignment for Swedish Conditions), funded by the Swedish Transport Administration (Trafikverket) under grant TRV 2023/109327.

\section*{Declaration of generative AI and AI-assisted technologies}

During the preparation of this work, the author used Claude (Anthropic) and ChatGPT (OpenAI) to assist with model implementation, code development and testing, and drafting and editing the manuscript.
All numerical results were produced by the model code.
The author reviewed and edited all content and takes full responsibility for the publication.

\newpage
\printbibliography[title={References}]

@article{abate2019,
  author  = {Abate, Megersa and Vierth, Inge and Karlsson, Rune and de Jong, Gerard and Baak, Jaap},
  title   = {A Disaggregate Stochastic Freight Transport Model for Sweden},
  journal = {Transportation},
  year    = {2019},
  volume  = {46},
  number  = {3},
  pages   = {671--696},
  doi     = {10.1007/s11116-018-9856-9}
}

@book{anderson1992,
  author    = {Anderson, Simon P. and de Palma, Andr{\'e} and Thisse, Jacques-Fran{\c{c}}ois},
  title     = {Discrete Choice Theory of Product Differentiation},
  publisher = {MIT Press},
  year      = {1992},
  address   = {Cambridge, MA},
  doi       = {10.7551/mitpress/2450.001.0001}
}

@article{baumol1970,
  author  = {Baumol, William J. and Bradford, David F.},
  title   = {Optimal Departures from Marginal Cost Pricing},
  journal = {American Economic Review},
  year    = {1970},
  volume  = {60},
  number  = {3},
  pages   = {265--283},
  url     = {https://www.jstor.org/stable/1817977}
}

@article{brouer2014,
  author  = {Brouer, Berit Dangaard and Alvarez, J. Fernando and Plum, Christian E. M. and Pisinger, David and Sigurd, Mikkel M.},
  title   = {A Base Integer Programming Model and Benchmark Suite for Liner-Shipping Network Design},
  journal = {Transportation Science},
  year    = {2014},
  volume  = {48},
  number  = {2},
  pages   = {281--312},
  doi     = {10.1287/trsc.2013.0471}
}

@article{dejong2007,
  author  = {de Jong, Gerard and Ben-Akiva, Moshe},
  title   = {A Micro-Simulation Model of Shipment Size and Transport Chain Choice},
  journal = {Transportation Research Part B: Methodological},
  year    = {2007},
  volume  = {41},
  number  = {9},
  pages   = {950--965},
  doi     = {10.1016/j.trb.2007.05.002}
}

@article{koza2020,
  author  = {Koza, David F. and Desaulniers, Guy and Ropke, Stefan},
  title   = {Integrated Liner Shipping Network Design and Scheduling},
  journal = {Transportation Science},
  year    = {2020},
  volume  = {54},
  number  = {2},
  pages   = {512--533},
  doi     = {10.1287/trsc.2018.0888}
}

@techreport{lindgren2019,
  author      = {Lindgren, Samuel and Johansson, Magnus and Vierth, Inge and Karlsson, Rune and Baak, Jaap and de Jong, Gerard},
  title       = {A Stochastic Logistics Module in Samgods: Work Tasks for Full Implementation},
  institution = {Swedish National Road and Transport Research Institute (VTI)},
  type        = {VTI PM},
  year        = {2019},
  month       = sep,
  note        = {Diarienummer 2018/0370-7.4},
  url         = {https://urn.kb.se/resolve?urn=urn:nbn:se:trafikverket:diva-11583}
}

@article{mankiw1986,
  author  = {Mankiw, N. Gregory and Whinston, Michael D.},
  title   = {Free Entry and Social Inefficiency},
  journal = {RAND Journal of Economics},
  year    = {1986},
  volume  = {17},
  number  = {1},
  pages   = {48--58},
  doi     = {10.2307/2555627},
  url     = {https://www.jstor.org/stable/2555627}
}

@incollection{mcfadden1974,
  author    = {McFadden, Daniel},
  title     = {Conditional Logit Analysis of Qualitative Choice Behavior},
  editor    = {Zarembka, Paul},
  booktitle = {Frontiers in Econometrics},
  publisher = {Academic Press},
  address   = {New York},
  year      = {1974},
  pages     = {105--142},
  url       = {https://eml.berkeley.edu/reprints/mcfadden/zarembka.pdf}
}

@article{morrow2011,
  author  = {Morrow, W. Ross and Skerlos, Steven J.},
  title   = {Fixed-Point Approaches to Computing Bertrand--Nash Equilibrium Prices under Mixed-Logit Demand},
  journal = {Operations Research},
  year    = {2011},
  volume  = {59},
  number  = {2},
  pages   = {328--345},
  doi     = {10.1287/opre.1100.0894}
}

@article{ramsey1927,
  author  = {Ramsey, Frank P.},
  title   = {A Contribution to the Theory of Taxation},
  journal = {Economic Journal},
  year    = {1927},
  volume  = {37},
  number  = {145},
  pages   = {47--61},
  doi     = {10.2307/2222721}
}

@article{tamannaei2021,
  author  = {Tamannaei, Mohammad and Zarei, Hamid and Aminzadegan, Sajede},
  title   = {A Game-Theoretic Approach to the Freight Transportation Pricing Problem in the Presence of Intermodal Service Providers in a Competitive Market},
  journal = {Networks and Spatial Economics},
  year    = {2021},
  volume  = {21},
  number  = {1},
  pages   = {123--173},
  doi     = {10.1007/s11067-020-09511-8}
}

@book{train2009,
  author    = {Train, Kenneth E.},
  title     = {Discrete Choice Methods with Simulation},
  edition   = {2},
  publisher = {Cambridge University Press},
  year      = {2009},
  doi       = {10.1017/CBO9780511805271}
}

@article{wang2014,
  author  = {Wang, Hua and Meng, Qiang and Zhang, Xiaoning},
  title   = {Game-Theoretical Models for Competition Analysis in a New Emerging Liner Container Shipping Market},
  journal = {Transportation Research Part B: Methodological},
  year    = {2014},
  volume  = {70},
  pages   = {201--227},
  doi     = {10.1016/j.trb.2014.09.006}
}

@article{chow2010,
  author  = {Chow, Joseph Y. J. and Yang, Choon Heon and Regan, Amelia C.},
  title   = {State-of-the-Art of Freight Forecast Modeling: Lessons Learned and the Road Ahead},
  journal = {Transportation},
  year    = {2010},
  volume  = {37},
  number  = {6},
  pages   = {1011--1030},
  doi     = {10.1007/s11116-010-9281-1}
}

@article{dejong2016,
  author  = {de Jong, Gerard and Tavasszy, L{\'o}r{\'a}nt and Bates, John and Gr{\o}nland, Stein Erik and Huber, Stefan and Kleven, Oskar and Lange, Peter and Ottem{\"o}ller, Ole and Schmorak, Nora},
  title   = {The Issues in Modelling Freight Transport at the National Level},
  journal = {Case Studies on Transport Policy},
  year    = {2016},
  volume  = {4},
  number  = {1},
  pages   = {13--21},
  doi     = {10.1016/j.cstp.2015.08.002}
}

@article{jensen2019,
  author  = {Jensen, Anders Fjendbo and Thorhauge, Mikkel and de Jong, Gerard and Rich, Jeppe and Dekker, Thijs and Johnson, Daniel and Ojeda Cabral, Manuel and Bates, John and Nielsen, Otto Anker},
  title   = {A Disaggregate Freight Transport Chain Choice Model for Europe},
  journal = {Transportation Research Part E: Logistics and Transportation Review},
  year    = {2019},
  volume  = {121},
  pages   = {43--62},
  doi     = {10.1016/j.tre.2018.10.004}
}

@article{keya2019,
  author  = {Keya, Nowreen and Anowar, Sabreena and Eluru, Naveen},
  title   = {Joint Model of Freight Mode Choice and Shipment Size: A Copula-Based Random Regret Minimization Framework},
  journal = {Transportation Research Part E: Logistics and Transportation Review},
  year    = {2019},
  volume  = {125},
  pages   = {97--115},
  doi     = {10.1016/j.tre.2019.03.007}
}

@article{berry1995,
  author  = {Berry, Steven and Levinsohn, James and Pakes, Ariel},
  title   = {Automobile Prices in Market Equilibrium},
  journal = {Econometrica},
  year    = {1995},
  volume  = {63},
  number  = {4},
  pages   = {841--890},
  doi     = {10.2307/2171802}
}

@article{aksoypierson2013,
  author  = {Aksoy-Pierson, Margaret and Allon, Gad and Federgruen, Awi},
  title   = {Price Competition Under Mixed Multinomial Logit Demand Functions},
  journal = {Management Science},
  year    = {2013},
  volume  = {59},
  number  = {8},
  pages   = {1817--1835},
  doi     = {10.1287/mnsc.1120.1664}
}

@article{gallego2014,
  author  = {Gallego, Guillermo and Wang, Ruxian},
  title   = {Multiproduct Price Optimization and Competition Under the Nested Logit Model with Product-Differentiated Price Sensitivities},
  journal = {Operations Research},
  year    = {2014},
  volume  = {62},
  number  = {2},
  pages   = {450--461},
  doi     = {10.1287/opre.2013.1249}
}

@article{gilbert2015,
  author  = {Gilbert, Fran{\c{c}}ois and Marcotte, Patrice and Savard, Gilles},
  title   = {A Numerical Study of the Logit Network Pricing Problem},
  journal = {Transportation Science},
  year    = {2015},
  volume  = {49},
  number  = {3},
  pages   = {706--719},
  doi     = {10.1287/trsc.2014.0560}
}

@article{shah2012,
  author  = {Shah, Nilopa and Brueckner, Jan K.},
  title   = {Price and Frequency Competition in Freight Transportation},
  journal = {Transportation Research Part A: Policy and Practice},
  year    = {2012},
  volume  = {46},
  number  = {6},
  pages   = {938--953},
  doi     = {10.1016/j.tra.2012.02.014}
}

@article{song2023,
  author  = {Song, Wei and Yang, Hualong and Li, Dechang and Yang, Zaili},
  title   = {Dynamic Pricing and Competition of Container Liner Shipping Services in a Duopoly Spot Market},
  journal = {Computers \& Industrial Engineering},
  year    = {2023},
  volume  = {185},
  pages   = {109613},
  doi     = {10.1016/j.cie.2023.109613}
}

@article{lyu2023,
  author  = {Lyu, Guangmei and Zhao, Min and Ji, Qingkai and Lin, Xiangmei},
  title   = {Improving Resilience via Capacity Allocation and Strategic Pricing: Co-opetition in a Shipping Supply Chain},
  journal = {Ocean \& Coastal Management},
  year    = {2023},
  volume  = {244},
  pages   = {106779},
  doi     = {10.1016/j.ocecoaman.2023.106779}
}

@article{wang2024alliance,
  author  = {Wang, Junjin and Xu, Songjun and Wu, Lingxiao and Wu, Sihong and Liu, Yanping},
  title   = {Multidimensional Container Shipping Alliance Decisions Among Competitors: Impact of Capacity Constraints and Market Competition},
  journal = {Transportation Research Part E: Logistics and Transportation Review},
  year    = {2024},
  volume  = {190},
  pages   = {103694},
  doi     = {10.1016/j.tre.2024.103694}
}

@article{lu2025,
  author  = {Lu, Bo and Fan, Lijie and Perera, Sandun C. and Wang, Jian-Jun},
  title   = {Does Antitrust Immunity Matter for Complementary Shipping Alliances? Competition and Welfare Analysis},
  journal = {Transportation Research Part E: Logistics and Transportation Review},
  year    = {2025},
  volume  = {198},
  pages   = {104125},
  doi     = {10.1016/j.tre.2025.104125}
}

@article{agarwal2008,
  author  = {Agarwal, Richa and Ergun, {\"O}zlem},
  title   = {Ship Scheduling and Network Design for Cargo Routing in Liner Shipping},
  journal = {Transportation Science},
  year    = {2008},
  volume  = {42},
  number  = {2},
  pages   = {175--196},
  doi     = {10.1287/trsc.1070.0205}
}

@article{brouer2014matheuristic,
  author  = {Brouer, Berit Dangaard and Desaulniers, Guy and Pisinger, David},
  title   = {A Matheuristic for the Liner Shipping Network Design Problem},
  journal = {Transportation Research Part E: Logistics and Transportation Review},
  year    = {2014},
  volume  = {72},
  pages   = {42--59},
  doi     = {10.1016/j.tre.2014.09.012}
}

@article{meng2014,
  author  = {Meng, Qiang and Wang, Shuaian and Andersson, Henrik and Thun, Kristian},
  title   = {Containership Routing and Scheduling in Liner Shipping: Overview and Future Research Directions},
  journal = {Transportation Science},
  year    = {2014},
  volume  = {48},
  number  = {2},
  pages   = {265--280},
  doi     = {10.1287/trsc.2013.0461}
}

@article{karsten2017,
  author  = {Karsten, Christian Vad and Brouer, Berit Dangaard and Pisinger, David},
  title   = {Competitive Liner Shipping Network Design},
  journal = {Computers \& Operations Research},
  year    = {2017},
  volume  = {87},
  pages   = {125--136},
  doi     = {10.1016/j.cor.2017.05.018}
}

@article{christiansen2020,
  author  = {Christiansen, Marielle and Hellsten, Erik Orm and Pisinger, David and Sacramento, David and Vilhelmsen, Charlotte},
  title   = {Liner Shipping Network Design},
  journal = {European Journal of Operational Research},
  year    = {2020},
  volume  = {286},
  number  = {1},
  pages   = {1--20},
  doi     = {10.1016/j.ejor.2019.09.057}
}

@article{lederer1993,
  author  = {Lederer, Phillip J.},
  title   = {A Competitive Network Design Problem with Pricing},
  journal = {Transportation Science},
  year    = {1993},
  volume  = {27},
  number  = {1},
  pages   = {25--38},
  doi     = {10.1287/trsc.27.1.25}
}

@article{gelareh2010,
  author  = {Gelareh, Shahin and Nickel, Stefan and Pisinger, David},
  title   = {Liner Shipping Hub Network Design in a Competitive Environment},
  journal = {Transportation Research Part E: Logistics and Transportation Review},
  year    = {2010},
  volume  = {46},
  number  = {6},
  pages   = {991--1004},
  doi     = {10.1016/j.tre.2010.05.005}
}

@article{angeloudis2016,
  author  = {Angeloudis, Panagiotis and Greco, Luciano and Bell, Michael G. H.},
  title   = {Strategic Maritime Container Service Design in Oligopolistic Markets},
  journal = {Transportation Research Part B: Methodological},
  year    = {2016},
  volume  = {90},
  pages   = {22--37},
  doi     = {10.1016/j.trb.2016.04.010}
}

@article{tawfik2019,
  author  = {Tawfik, Christine and Limbourg, Sabine},
  title   = {A Bilevel Model for Network Design and Pricing Based on a Level-of-Service Assessment},
  journal = {Transportation Science},
  year    = {2019},
  volume  = {53},
  number  = {6},
  pages   = {1609--1626},
  doi     = {10.1287/trsc.2019.0906}
}

@article{wang2023integrated,
  author  = {Wang, Zhenjie and Zhang, Dezhi and Tavasszy, L{\'o}r{\'a}nt and Fazi, Stefano},
  title   = {Integrated Multimodal Freight Service Network Design and Pricing with a Competing Service Integrator and Heterogeneous Shipper Classes},
  journal = {Transportation Research Part E: Logistics and Transportation Review},
  year    = {2023},
  volume  = {179},
  pages   = {103290},
  doi     = {10.1016/j.tre.2023.103290}
}

@article{bresnahan1991,
  author  = {Bresnahan, Timothy F. and Reiss, Peter C.},
  title   = {Entry and Competition in Concentrated Markets},
  journal = {Journal of Political Economy},
  year    = {1991},
  volume  = {99},
  number  = {5},
  pages   = {977--1009},
  doi     = {10.1086/261786}
}

@article{berry1992,
  author  = {Berry, Steven T.},
  title   = {Estimation of a Model of Entry in the Airline Industry},
  journal = {Econometrica},
  year    = {1992},
  volume  = {60},
  number  = {4},
  pages   = {889--917},
  doi     = {10.2307/2951571}
}

@article{ericson1995,
  author  = {Ericson, Richard and Pakes, Ariel},
  title   = {Markov-Perfect Industry Dynamics: A Framework for Empirical Work},
  journal = {Review of Economic Studies},
  year    = {1995},
  volume  = {62},
  number  = {1},
  pages   = {53--82},
  doi     = {10.2307/2297841}
}

@article{ciliberto2009,
  author  = {Ciliberto, Federico and Tamer, Elie},
  title   = {Market Structure and Multiple Equilibria in Airline Markets},
  journal = {Econometrica},
  year    = {2009},
  volume  = {77},
  number  = {6},
  pages   = {1791--1828},
  doi     = {10.3982/ECTA5368}
}

@article{aguirregabiria2012,
  author  = {Aguirregabiria, Victor and Ho, Chun-Yu},
  title   = {A Dynamic Oligopoly Game of the US Airline Industry: Estimation and Policy Experiments},
  journal = {Journal of Econometrics},
  year    = {2012},
  volume  = {168},
  number  = {1},
  pages   = {156--173},
  doi     = {10.1016/j.jeconom.2011.09.013}
}

@article{horcher2021,
  author  = {H{\"o}rcher, Daniel and Tirachini, Alejandro},
  title   = {A Review of Public Transport Economics},
  journal = {Economics of Transportation},
  year    = {2021},
  volume  = {25},
  pages   = {100196},
  doi     = {10.1016/j.ecotra.2021.100196}
}

@techreport{asplundwestin2017,
  author      = {Asplund, Disa and Westin, Jonas},
  title       = {Modellering av slingor inom sjötransporter},
  institution = {Statens väg- och transportforskningsinstitut (VTI)},
  type        = {VTI notat},
  number      = {31-2016},
  year        = {2017},
  langid      = {swedish},
  url         = {https://urn.kb.se/resolve?urn=urn:nbn:se:trafikverket:diva-6011}
}

@techreport{johansson2021,
  author      = {Johansson, Magnus and Karlsson, Rune and Westin, Jonas},
  title       = {Sjöfart mot svensk hamn 2007 till 2019: förändringar i volym, struktur och rörelsemönster},
  institution = {Statens väg- och transportforskningsinstitut (VTI)},
  type        = {VTI rapport},
  number      = {1070},
  year        = {2021},
  langid      = {swedish},
  url         = {https://urn.kb.se/resolve?urn=urn:nbn:se:trafikverket:diva-5786}
}

@inproceedings{westin2024loops,
  author    = {Westin, Jonas and Olsson, Leif and Åhag, Per},
  title     = {Loop-Based Maritime Transport Optimization in Swedish National Freight Transport Models},
  booktitle = {2024 IEEE International Conference on Industrial Engineering and Engineering Management (IEEM)},
  year      = {2024},
  pages     = {1--5},
  publisher = {IEEE},
  doi       = {10.1109/IEEM62345.2024.10857080}
}

@article{westin2026loops,
  author  = {Westin, Jonas and Åhag, Per and Olsson, Leif and Asplund, Disa},
  title   = {Modeling Loops and Backhaul in Maritime Freight: A MILP Extension for Large-Scale Freight Transport Models},
  journal = {Transportation Research Interdisciplinary Perspectives},
  volume  = {38},
  pages   = {102059},
  year    = {2026},
  doi     = {10.1016/j.trip.2026.102059}
}

@incollection{benakivadejong2013,
  author    = {Ben-Akiva, Moshe and de Jong, Gerard},
  title     = {The Aggregate--Disaggregate--Aggregate (ADA) Freight Model System},
  booktitle = {Freight Transport Modelling},
  editor    = {Ben-Akiva, Moshe and Meersman, Hilde and Van de Voorde, Eddy},
  year      = {2013},
  pages     = {69--90},
  publisher = {Emerald Group Publishing},
  doi       = {10.1108/9781781902868-004}
}

@techreport{dejongbaak2023,
  author      = {de Jong, Gerard and Baak, Jaap},
  title       = {Method Report -- Logistics Model in the Swedish National Freight Model System},
  institution = {Significance},
  type        = {Deliverable for Trafikverket},
  year        = {2024},
  month       = jan,
  url         = {https://bransch.trafikverket.se/contentassets/ab220f9016154ef7a8478555560bb280/2024/method-report---logistics-model-in-the-swedish-national-freight-model-system.pdf}
}

@techreport{trafikverket2026,
  author      = {{Trafikverket}},
  title       = {Samgods 1.2.3 -- User Manual},
  institution = {Trafikverket},
  type        = {Report},
  year        = {2026},
  month       = may,
  note        = {Prepared by G. Sala (SYSTRA) and the Samgods working group at Trafikverket},
  url         = {https://bransch.trafikverket.se/contentassets/ab220f9016154ef7a8478555560bb280/2026/samgods-1.2.3---user-manual.pdf},
  langid      = {english}
}

@article{caplin1991,
  author  = {Caplin, Andrew and Nalebuff, Barry},
  title   = {Aggregation and Imperfect Competition: On the Existence of Equilibrium},
  journal = {Econometrica},
  year    = {1991},
  volume  = {59},
  number  = {1},
  pages   = {25--59},
  doi     = {10.2307/2938239}
}

@article{facchinei2007,
  author  = {Facchinei, Francisco and Kanzow, Christian},
  title   = {Generalized {Nash} Equilibrium Problems},
  journal = {4OR},
  year    = {2007},
  volume  = {5},
  number  = {3},
  pages   = {173--210},
  doi     = {10.1007/s10288-007-0054-4}
}

@article{lihuh2011,
  author  = {Li, Hongmin and Huh, Woonghee Tim},
  title   = {Pricing Multiple Products with the Multinomial Logit and Nested Logit Models: Concavity and Implications},
  journal = {Manufacturing \& Service Operations Management},
  year    = {2011},
  volume  = {13},
  number  = {4},
  pages   = {549--563},
  doi     = {10.1287/msom.1110.0344}
}

@article{everett1963,
  author  = {Everett, Hugh},
  title   = {Generalized {Lagrange} Multiplier Method for Solving Problems of Optimum Allocation of Resources},
  journal = {Operations Research},
  year    = {1963},
  volume  = {11},
  number  = {3},
  pages   = {399--417},
  doi     = {10.1287/opre.11.3.399}
}

\clearpage
\appendix

\section{Parameters}
\label{app:parameters}

\Cref{tab:parameters} lists the inputs of the six-port instance, generated from the input workbook and the settings of the logged run. Ranges are given where a parameter varies by port, leg or market; the full OD tables are in the input workbook distributed with the code. The nested-logit runs of \cref{sec:nested} change only the nest parameter.

\begin{table}[H]
  \centering
  \caption{Parameters of the six-port test instance.}
  \label{tab:parameters}
  \small
\begin{tabularx}{\textwidth}{>{\raggedright\arraybackslash}p{0.27\textwidth} l X l}
    \toprule
    Parameter & Symbol & Value & Unit \\
    \midrule
    Ports & -- & STO, HEL, TAL, RIX, LBC, RTM & -- \\
    Maritime distances & $d_{ij}$ & 45--800 & nautical miles \\
    Port time per call & $h_i$ & 6--10 & hours \\
    Port-call cost & $C^{\mathrm{port}}_i$ & 12,000--26,000 & per call \\
    Vessel capacity by class & $u$ & SMALL: 500, LARGE: 900 & FFE \\
    Weekly vessel cost by class & $C^{\mathrm{charter}}_f$ & SMALL: 120,000, LARGE: 185,000 & per week \\
    Fuel coefficient by class & $a_f$ & SMALL: 0.003, LARGE: 0.0045 & t/day/kn$^3$ \\
    Speed levels & $v$ & 14, 17, 20 & knots \\
    Fuel price & $q^{\mathrm{fuel}}$ & 900 & per tonne \\
    Variable cost & $c^{\mathrm{var}}_{kg}$ & HIGH: 70, LOW: 50 & per FFE \\
    Weekly demand & $D_{kg}$ & 50--190 (total 2,295) & FFE per week \\
    Alternative-specific constant & $\gamma_{0,g}$ & HIGH: 4.8, LOW: 4.2 & utils \\
    Transit-time coefficient & $\gamma_{T,g}$ & HIGH: 0.045, LOW: 0.018 & utils per hour \\
    Price coefficient & $\beta_g$ & HIGH: 0.0018, LOW: 0.0024 & utils per currency unit \\
    Nest parameter (sea nest) & $\rho_g$ & 1 & -- \\
    Outside-option utility & $V^{\mathrm{out}}_{kg}$ & 0.80--2.15 & utils \\
    Maximum transit time & $TT_k$ & 24--96 & hours \\
    Price bounds & -- & [-5,000, 25,000] & per FFE \\
    Utility normalisation & -- & utilities as in the input; $CS$ is the logsum (eq.~\ref{eq:consumer-surplus-market}) & -- \\
    Maximum port calls per rotation & -- & 4 & -- \\
    Tolerances & $\varepsilon_p$, $\varepsilon_{\pi,p}$, $\varepsilon_\pi$, $\varepsilon_J$ & 1e-05, 0.0001, 1, 1 & -- \\
    Damping & $\theta$ & 1 (fallback 0.5) & -- \\
    Best-response grid & -- & 160 points (certificate 1200), geometric in the mark-up from 0.05 & -- \\
    \bottomrule
\end{tabularx}
 \end{table}

\newpage
\section{Competitive solution algorithm}
\label{app:algorithm}

\Cref{alg:competition} summarises the solve--verify--repair procedure of \cref{sec:solution} for one initial service set. Replacements and entries are evaluated in decreasing order of their estimated gain, and a candidate $r$ is skipped when the bound \eqref{eq:standalone-bound} shows that it cannot be profitable, that is, when $B_r\le\varepsilon_\pi$ for entry or $B_r-\pi_s^*(A)\le\varepsilon_\pi$ for replacing $s$. Clone entries are evaluated during verification as a diagnostic. The procedure is run from every starting configuration, and the verified endpoints are then grouped into distinct equilibria.

\begin{algorithm}[htbp]
\caption{Solve--verify--repair procedure for the competitive regime}
\label{alg:competition}
\begin{algorithmic}[1]
\Require candidate set $S$, initial active set $A_0\subseteq S$, tolerance $\varepsilon_\pi$, bounds $B_r$ for $r\in S$
\Ensure verified equilibrium $(A^*,p^*)$, or a reported cycle
\State $A\gets A_0$, \ $\mathcal V\gets\emptyset$ \Comment{visited sets}
\Loop
  \State solve $p^*(A)$ by Gauss--Seidel global best responses from the stand-alone prices
  \State certify $p^*(A)$ on the fine grid
  \If{$\min_{s\in A}\pi_s^*(A)<-\varepsilon_\pi$}
    \State $A'\gets A\setminus\{\arg\min_{s\in A}\pi_s^*(A)\}$ \Comment{exit}
  \ElsIf{a replacement $s\to r$ with gain $>\varepsilon_\pi$ is found}
    \State $A'\gets$ the set after the best replacement \Comment{replacement}
  \ElsIf{an entrant $r$ with $\pi_r^*(A\cup\{r\})>\varepsilon_\pi$ is found}
    \State $A'\gets A\cup\{\text{best entrant}\}$ \Comment{entry}
  \Else
    \State verify all exits, entries and replacements not excluded by $B_r$
    \If{a deviation with gain $>\varepsilon_\pi$ exists}
      \State $A'\gets$ the set after the most profitable deviation \Comment{repair}
    \Else
      \State \Return $(A,p^*(A))$ as a verified equilibrium
    \EndIf
  \EndIf
  \If{$A'\in\mathcal V$}
    \State \Return cycle through $A'$
  \EndIf
  \State $\mathcal V\gets\mathcal V\cup\{A\}$, \ $A\gets A'$
\EndLoop
\end{algorithmic}
\end{algorithm}

\section{Verification diagnostics}
\label{app:diagnostics}

\Cref{tab:diagnostics} reports, for each scenario, the size of the verification problem and the quality of the reported competitive equilibria under multinomial logit, and \cref{tab:nest-diagnostics} the same for the nested-logit run. `Pruned by bound' counts deviations excluded by \eqref{eq:standalone-bound}; every other deviation was solved. Run times depend on the machine.

\begin{table}[H]
  \centering
  \caption{Verification diagnostics for the competitive regime (multinomial logit).}
  \label{tab:diagnostics}
  \small
\begin{tabularx}{\textwidth}{>{\raggedright\arraybackslash}X r r r r }
    \toprule
    \textbf{Diagnostic} & \textbf{Base} & \textbf{High fuel} & \textbf{High demand} & \textbf{Low capacity} \\
    \midrule
    Candidate services $|S|$ (enumerated) & 758 (870) & 758 (870) & 758 (870) & 758 (870) \\
    Services that can break even & 201 & 154 & 318 & 201 \\
    Entry deviations tested / pruned by bound & 197 / 557 & 151 / 604 & 313 / 440 & 197 / 557 \\
    Replacement deviations tested / pruned & 487 / 2529 & 208 / 2057 & 1031 / 2734 & 487 / 2529 \\
    Clone deviations tested & 4 & 3 & 5 & 4 \\
    Largest clone-entry profit & \num{-38943} & \num{-42036} & \num{-16872} & \num{-38943} \\
    Largest unilateral gain & \num{-273.4} & \num{-851.3} & \num{-433.3} & \num{-273.4} \\
    Largest price residual $r_p$ & \num{7.8e-07} & \num{8.4e-07} & \num{1.6e-07} & \num{7.8e-07} \\
    Largest price-stage gain $r_\pi$ & \num{0.0e+00} & \num{0.0e+00} & \num{5.8e-11} & \num{0.0e+00} \\
    Largest leg load factor & 0.476 & 0.516 & 0.546 & 0.953 \\
    Failed price solves & 0 & 0 & 0 & 0 \\
    Starts / distinct verified equilibria & 110 / 1 & 110 / 2 & 110 / 1 & 110 / 1 \\
    Cycles detected & 0 & 0 & 0 & 0 \\
    Price equilibria found for reported set & 1 & 1 & 1 & 1 \\
    Run time (s) & 129 & 62 & 286 & 207 \\
    \bottomrule
\end{tabularx}
 \end{table}

\begin{table}[H]
  \centering
  \caption{Verification diagnostics for the competitive regime (nested logit, $\rho=0.6$).}
  \label{tab:nest-diagnostics}
  \small
\begin{tabularx}{\textwidth}{>{\raggedright\arraybackslash}X r r r r }
    \toprule
    \textbf{Diagnostic} & \textbf{Base} & \textbf{High fuel} & \textbf{High demand} & \textbf{Low capacity} \\
    \midrule
    Candidate services $|S|$ (enumerated) & 758 (870) & 758 (870) & 758 (870) & 758 (870) \\
    Services that can break even & 201 & 154 & 318 & 201 \\
    Entry deviations tested / pruned by bound & 198 / 557 & 151 / 604 & 315 / 440 & 198 / 557 \\
    Replacement deviations tested / pruned & 184 / 2081 & 189 / 2076 & 194 / 2071 & 184 / 2081 \\
    Clone deviations tested & 3 & 3 & 3 & 3 \\
    Largest clone-entry profit & \num{-55019} & \num{-57885} & \num{-24398} & \num{-55019} \\
    Largest unilateral gain & \num{-1138.1} & \num{-851.3} & \num{-2386.0} & \num{-1138.1} \\
    Largest price residual $r_p$ & \num{2.3e-13} & \num{8.4e-07} & \num{2.3e-13} & \num{2.3e-13} \\
    Largest price-stage gain $r_\pi$ & \num{5.8e-11} & \num{2.9e-11} & \num{5.8e-11} & \num{5.8e-11} \\
    Largest leg load factor & 0.516 & 0.516 & 0.645 & 1.000 \\
    Failed price solves & 0 & 0 & 0 & 0 \\
    Starts / distinct verified equilibria & 110 / 4 & 110 / 4 & 110 / 8 & 110 / 4 \\
    Cycles detected & 0 & 0 & 0 & 0 \\
    Price equilibria found for reported set & 1 & 1 & 1 & 1 \\
    Run time (s) & 83 & 64 & 174 & 212 \\
    \bottomrule
\end{tabularx}
 \end{table}

\end{document}